\documentclass[aop,MSNbibl,citesort,seceqn,dvips]{arximspdf}
\usepackage{graphicx}

\doi{10.1214/10-AOP624}
\volume{40}
\issue{2}
\pubyear{2012}
\firstpage{535}
\lastpage{577}

\makeatletter

\newtheorem{theorem}{Theorem}[section]
\newtheorem{lemma}[theorem]{Lemma}
\newtheorem{proposition}[theorem]{Proposition}

\newproclaim{assumption}[theorem]{Assumption}
\newproclaim{remark}[theorem]{Remark}

\newcommand{\iint}{\int\!\!\int}

\newcommand{\E}{\mathbf{E}}
\newcommand{\bL}{\mathbf{L}}
\newcommand{\N}{\mathbf{N}}
\newcommand{\Z}{\mathbf{Z}}
\newcommand{\p}{\mathbf{P}}

\newcommand{\CA}{\mathcal{A}}
\newcommand{\CC}{\mathcal{C}}
\newcommand{\CL}{\mathcal{L}}
\newcommand{\CM}{\mathcal{M}}
\newcommand{\CR}{\mathcal{R}}
\newcommand{\CS}{\mathcal{S}}
\newcommand{\CT}{\mathcal{T}}
\newcommand{\CX}{\mathcal{X}}

\newcommand{\CG}{\mathcal{G}}

\newcommand{\diam}{{\mathrm{diam}}}
\newcommand{\cov}{{\mathrm{cov}}}

\newcommand{\one}{\mathbf{1}}

\newcommand{\ol}{\overline}
\newcommand{\ul}{\underline}
\newcommand{\wt}{\widetilde}

\newcommand{\hit}{{\mathrm{hit}}}
\newcommand{\mix}{{\mathrm{mix}}}

\makeatother

\begin{document}
\begin{frontmatter}

\title{Uniformity of the uncovered set of random walk and cutoff for lamplighter chains}
\runtitle{Uniformity of the uncovered set of random walk}

\begin{aug}
\author[A]{\fnms{Jason} \snm{Miller}\corref{}\ead[label=e2]{jmiller@math.stanford.edu}} and
\author[B]{\fnms{Yuval} \snm{Peres}\ead[label=e3]{peres@microsoft.com}}
\runauthor{J. Miller and Y. Peres}
\affiliation{Stanford University and Microsoft Research}
\address[A]{Department of Mathematics\\
Stanford University\\
Stanford, California 94305\\
USA\\
\printead{e2}}
\address[B]{Microsoft Research\\
Theory Group\\
Redmond, Washington 98052\\
USA\\
\printead{e3}}
\end{aug}

\received{\smonth{1} \syear{2010}}
\revised{\smonth{10} \syear{2010}}

%
\begin{abstract}
We show that the measure on markings of $\Z_n^d$, $d \geq3$, with
elements of $\{0,1\}$ given by i.i.d. fair coin flips on the range
$\CR$ of a~random walk $X$ run until time $T$ and $0$ otherwise becomes
indistinguishable from the uniform measure on such markings at the
threshold $T = \frac{1}{2}T_\cov(\Z_n^d)$. As a consequence of our
methods, we show that the total variation mixing time of the random
walk on the lamplighter graph $\Z_2 \wr\Z_n^d$, $d \geq3$, has a
cutoff with threshold $\frac{1}{2} T_\cov(\Z_n^d)$. We give a~general
criterion under which both of these results hold; other examples for
which this applies include bounded degree expander families, the
intersection of an infinite supercritical percolation cluster with an
increasing family of balls, the hypercube and the Caley graph of the
symmetric group generated by transpositions. The proof also yields
precise asymptotics for the decay of correlation in the uncovered set.
\end{abstract}

%

\end{frontmatter}

\section{Introduction}\label{sec1}

\begin{figure}
\begin{tabular}{@{}c@{\hspace*{5pt}}c@{\hspace*{5pt}}c@{}}

\includegraphics{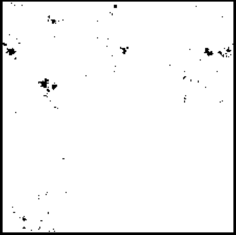}
 & \includegraphics{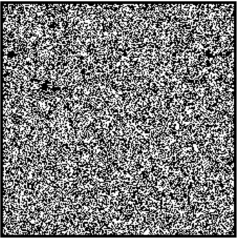} & \includegraphics{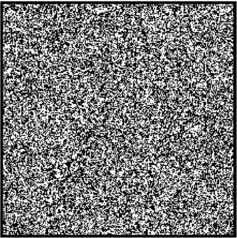}\\
(a)  & (b) & (c)
\end{tabular}
\caption{The subset $\CL(\frac{1}{2}, \Z_n^2)$ of $\Z_n^2$ consisting
of those points unvisited by a random walk~$X$ run for $\frac{1}{2}
T_\cov(\Z_n^2)$, where $T_\cov(\Z_n^2)$ is the expected number of steps
required for~$X$ to cover $\Z_n^2$, exhibits clustering. Consequently,
the marking of $\Z_n^2$ by elements of $\{0,1\}$ given by the results
of i.i.d. coin flips on the range of $X$ at time $\frac{1}{2}
T_\cov(\Z_n^2)$ and zero otherwise can be distinguished from a uniform
marking. \textup{(a)} $\CL(\frac{1}{2},\Z_n^2)$. \textup{(b)} $\CL(\frac{1}{2}, \Z_n^2)$
marked with i.i.d. coin flips. \textup{(c)} $\Z_n^2$ marked with i.i.d. coin
flips.}\label{fig1}
\end{figure}

Suppose $G = (V,E)$ is a finite, connected graph and $X$ is a lazy
random walk on $G$. This means that $X$ is the Markov chain with state
space $V$ and transition kernel
\[
p(x,y;G) = \p_x[X(1) = y] = \cases{
\displaystyle \frac{1}{2}, &\quad if $x = y$,\cr
\displaystyle \frac{1}{2\deg(x)}, &\quad if $\{x,y\} \in E$.}
\]
Let
\[
\tau_\cov(G) = \min\bigl\{t \geq0 \dvtx V \mbox{ is contained in the range
of } X|_{[0,t]}\bigr\}
\]
be the \textit{cover time} and let $T_\cov(G) = \E_\pi[\tau_\cov
(G)]$ be the \textit{expected cover time}. Here and hereafter, a
subscript of $\pi$ indicates that $X$ is started from stationarity.
Let $\tau(y) = \min\{t \geq0 \dvtx X(t) = y\}$ be the first time $X$
hits $y$ and
\[
T_\hit(G) = \max_{x,y \in V} \E_x[\tau(y)]
\]
be the \textit{maximal hitting time}. If $(G_n)$ is a sequence of graphs
with $T_\hit(G_n) = o(T_\cov(G_n))$, then a result of Aldous \cite
{A91}, Theorem 2, implies that $\tau_\cov(G_n)$ has a threshold around
its mean: $\tau_\cov(G_n) / T_\cov(G_n) = 1 + o(1)$. Many sequences
of graphs satisfy this condition, for example, $\Z_n^d$ for $d \geq2$,
$\Z_2^n$, and the complete graph $K_n$. When Aldous' condition holds,
the set
\[
\CL(\alpha;G_n) = \{x \in V_n \dvtx\tau(x) \geq\alpha T_\cov(G_n)\},
\]
$V_n$ the vertices of $G_n$, of \textit{$\alpha$-late points}, that is,
points hit after time $\alpha T_\cov(G_n)$, $\alpha\in(0,1)$, often
has an interesting structure. The case $G_n = \Z_n^2$ was first
studied by Brummelhuis and Hilhorst in \cite{BH91} where it is shown
that $\E|\CL(\alpha;\Z_n^2)|$ has growth exponent $2(1-\alpha)$
and that points in $\CL_n(\alpha;\Z_n^2)$ are positively correlated.
This suggests that $\CL(\alpha;G_n)$ has a fractal structure and
exhibits clustering. These statements were made precise by Dembo,
Peres, Rosen and Zeitouni in \cite{DPRZLATE06} where they show that
the growth exponent of $|\CL(\alpha;\Z_n^2)|$ is $2(1-\alpha)$
\textit{with high probability} in addition to making a rigorous
quantification of the clustering phenomenon
(see Figure \ref{fig1} for an illustration of this).

%

If $G_n$ is either $K_n$ or $\Z_n^d$ for $d \geq3$, then it is also
true that ${\log}|\CL(\alpha;G_n)| \sim{(1-\alpha) \log}|V_n|$ with
high probability. In contrast to $\CL(\alpha;\Z_n^2)$, $\CL(\alpha
;K_n)$ does not exhibit clustering and is ``uniformly random'' in the
sense that conditional on $s_0 = |\CL(\alpha;K_n)|$, all subsets of
$K_n$ of size $s_0$ are equally likely. The rapid decay of correlation
in $\CL(\alpha;\Z_n^d)$ for $d \geq3$ determined by Brummelhuis and
Hilhorst \cite{BH91} indicates that the clustering phenomenon is also
not present in this case and leads one to speculate that $\CL(\alpha
;\Z_n^d)$ is likewise in some sense ``uniformly random.''

The purpose of this article is to quantify the degree to which this
holds. We use as our measure of uniformity the following statistical
test. Let $\CR(\alpha;G)$ be the (random) subset of $V$ covered by
$X$ at time $\alpha T_\cov(G)$ and let $\mu(\cdot;\alpha,G)$ be the
probability measure on $\CX(G) = \{f \dvtx V \to\{0,1\}\}$ given by
first sampling~$\CR(\alpha;G)$ then setting
\[
f(x) = \cases{
\xi(x), &\quad if $x \in\CR(\alpha;G)$,\cr
0, &\quad otherwise,}
\]
where $(\xi(x) \dvtx x \in V)$ is a collection of i.i.d. variables such
that $\p[\xi(x) = 0] = \p[\xi(x) = 1] = \frac{1}{2}$. The
question we are interested in is:\vspace*{8pt}

\textit{How large does $\alpha\in(0,1)$ need to be so that $\mu(\cdot
;\alpha,G)$ is indistinguishable from the uniform measure $\nu(\cdot
;G)$ on $\CX(G)$}?\vspace*{8pt}

It must be that $\alpha\geq1/2$ in the case of $\Z_n^d$
for $d \geq2$ since if $\alpha< 1/2$ then
\[
\frac{|\CL(\alpha;\Z_n^d)| - ({1/2})n^d} {n^{d/2}} \to\infty
\qquad\mbox{as } n \to\infty.
\]
In particular, the deviations of the number of zeros from $n^d/2$ which
arise in a~marking from such $\alpha$ far exceed that in the uniform
case. By \cite{A91}, Theorem 2, it is also true that $\alpha\leq1$
since if $\alpha> 1$ then with high probability $|\CL(\alpha;\Z
_n^d)| = 0$. The main result of this article is that the threshold for
indistinguishability for any sequence of graphs $(G_n)$ with $\lim_{n
\to\infty} |V_n| = \infty$ is $\alpha= \frac{1}{2}$ provided
random walk on $(G_n)$ is \textit{uniformly locally transient} and
satisfies a mild connectivity hypothesis.

We need the following definitions in order to give a precise statement
of our results. The $\varepsilon$-\textit{total variation mixing time} of
$G$ is
\[
T_\mix(\varepsilon;G) = \min\Bigl\{t \geq0 \dvtx\max_{x \in V} \|
p^t(x,\cdot
;G) - \pi\|_{\mathrm{TV}} \leq\varepsilon\Bigr\},
\]
where $p^t(x,y;G) = \p_x[ X(t) = y]$ is the $t$-step transition kernel
of $X$ started at $x$,
\[
\| \mu- \nu\|_{\mathrm{TV}} = \max_{A \subseteq V} | \mu(A) - \nu(A)| =
\frac{1}{2} \sum_{x \in V} |\mu(x) - \nu(x)|
\]
is the \textit{total variation distance} between the measures $\mu,\nu
$ on $V$ and $\pi$ is the stationary distribution of $X$. The
$\varepsilon$-\textit{uniform mixing time} of $G$ is
\[
T_\mix^U(\varepsilon;G) = \min\biggl\{ t \geq0 \dvtx\max_{x,y \in V}
\biggl| \frac{p^t(x,y;G)}{\pi(y)} - 1 \biggr| \leq\varepsilon\biggr\}.
\]
It is a basic fact (\cite{AFMC,LPW08}; see also Proposition \ref
{propmixingdecay}) that $T_\mix^U(\varepsilon;G)$ is within a factor of
${\log}|V|$ of $T_\mix(\varepsilon;G)$, however, for many graphs this
factor is constant. Whenever we omit $\varepsilon$ and write $T_\mix(G),
T_\mix^U(G)$ it is understood that $\varepsilon= \frac{1}{4}$.
\textit{Green's function} of $G$ is
\[
g(x,y;G) = \sum_{t=0}^{T_\mix^U(G)} p^t(x,y;G),
\]
that is, the expected amount of time that $X$ spends at $y$ until time
$T_\mix^U(G)$ when started at $x$. For $A \subseteq V$, we set
\[
g(x,A;G) = \sum_{y \in A} g(x,y;G).
\]
We say that $(G_n)$ is \textit{uniformly locally transient} with \textit
{transience function}~$\rho\dvtx\allowbreak
[0$, $\infty) \times[0,\infty) \to
[0,\infty)$ if
\[
g(x,A;G_n) \leq\rho(d(x,A), \diam(A)) \qquad\mbox{for all } n \mbox{
and } x \in V_n, A \subseteq V_n.
\]
Here, $d(\cdot,\cdot)$ is the graph distance, $d(x,A) = \min_{y \in
A} d(x,y)$, and $\rho(\cdot,s)$ is assumed to be nonincreasing with
$\lim_{r \to\infty} \rho(r,s) = 0$ when $s$ is fixed. Let $\rho(r)
= \rho(r,1)$,
\[
\ol{\Delta}(G) = \max_{x \in V} \deg(x),\qquad \ul{\Delta}(G) = \min
_{x \in V} \deg(x)\quad \mbox{and}\quad \Delta(G) = \frac{\ol{\Delta
}(G)}{\ul{\Delta}(G)}.
\]
\begin{assumption}[(Transience)]
\label{assumpgraphs}
$(G_n)$ is a sequence of uniformly locally transient graphs with $|V_n|
\to\infty$ such that there exists $\Delta_0 > 0$ so that\break $\Delta
(G_n) \leq\Delta_0$ for all $n$ and, for each $r > 0$:
\begin{longlist}[(1)]
\item[(1)]\hypertarget{assumpgraphsball} 
${\log}|B(x,r)| = o({\log}|V_n|)$ as $n
\to\infty$, and\vspace*{2pt}
\item[(2)]\hypertarget{assumpgraphshit} $T_\mix^U(G_n) \ol{\Delta
}{}^r(G_n) =
o(|V_n|)$ as $n \to\infty$.
\end{longlist}
\end{assumption}

The reason for the hypothesis $\Delta(G_n) \leq\Delta_0$ is that it implies
\[
\frac{\pi(x;G_n)}{ \pi(y;G_n)} \leq\Delta_0 \qquad\mbox{uniformly in }
x,y \in V_n \mbox{ and } n.
\]
In particular, this combined with uniform local transience allows us to
conclude that the hitting time of any two points $x,y \in V_n$ is
comparable. The purpose of part \hyperlink{assumpgraphsball}{(1)} of Assumption
\ref{assumpgraphs} is to ensure that for every $r,n > 0$ we can
construct an $r$-net $E_{r,n}$ of $V_n$ whose size at logarithmic
scales is comparable to $|V_n|$, that is, ${\log}|E_{r,n}| = {\log}|V_n|
+ o(1)$ as $n \to\infty$. Finally, part \hyperlink
{assumpgraphshit}{(2)} of
Assumption \ref{assumpgraphs} is important since by a union bound it
implies that\vspace*{1pt} the probability that $X$ hits any fixed ball of finite
radius within time $T_\mix^U(G_n)$ when initialized from stationarity
tends to zero with $n$.

We will also need to make the following assumption.
\begin{assumption}[(Connectivity)]
\label{assumpfarawayharnack}
$(G_n)$ is a sequence of graphs satisfying either:
\begin{longlist}[(2)]
\item[(1)]\hypertarget{assumpfaraway}
for every $\gamma> 0$ there exists
$R_n^\gamma\to\infty$ as $n \to\infty$ satisfying $R_n^\gamma\leq
\frac{1}{2}\times\max\{ {R > 0 \dvtx\max_{x \in V_n}} |B(x,R)| \leq
|V_n|^{\gamma}\}$ such that for every $r > 0$,
\[
\frac{T_\mix^U(G_n)}{R_n^\gamma} \max_{d(x,A) \geq R_n^\gamma}
g(x,A) = o(1) \qquad\mbox{as } n \to\infty
\]
uniformly in $A \subseteq V_n$ with $\diam(A) \leq r$, or\vadjust{\goodbreak}
\item[(2)]
\hypertarget{assumpharnack} a uniform Harnack inequality, that is, for
each $\alpha> 1$ there exists $C = C(\alpha) > 0$ such that for every
$x,r,R > 0$ with $R/r \geq\alpha$ and positive harmonic function $h$
on $B(x,R)$ we have that
\[
\max_{y \in B(x,r)} h(y) \leq C \min_{y \in B(x,r)} h(y).
\]
\end{longlist}
\end{assumption}

Assumption \ref{assumpfarawayharnack} ensures that $(G_n)$ is in some
sense well connected. In particular, part \hyperlink
{assumpfaraway}{(1)} is
used to show that $X$ is uniformly unlikely to hit a~small ball before
remixing provided its starting point and the small ball are far enough
apart. This hypothesis will be relevant for graphs where $|\partial
B(x,r)|$ is comparable to or larger than $|B(x,r)|$, as in the case of
$\Z_2^n$ or graphs which are locally tree-like. Part \hyperlink
{assumpharnack}{(2)}
is meant to be applicable for graphs where $|\partial
B(x,r)|$ is much smaller than $|B(x,r)|$, as in the case of~$\Z_n^d$,
and is used to deduce that the empirical average of the probability
that successive excursions of $X$ between concentric spheres $\partial
B(x,r), \partial B(x,R)$ hit~$x$ conditional on their entrance and exit
points is well concentrated around its mean provided $R > r$ are large enough.

We now state our main theorem.
\begin{theorem}
\label{thmthreshold}
If $(G_n)$ satisfies Assumptions \ref{assumpgraphs} and \ref
{assumpfarawayharnack}, then for every \mbox{$\varepsilon> 0$},
\[
\lim_{n \to\infty} \biggl\| \mu\biggl(\cdot;\frac{1}{2} + \varepsilon,G_n\biggr) -
\nu(\cdot;G_n) \biggr\|_{\mathrm{TV}} = 0
\]
and
\[
\lim_{n \to\infty} \biggl\| \mu\biggl(\cdot;\frac{1}{2} - \varepsilon,G_n\biggr) -
\nu(\cdot;G_n) \biggr\|_{\mathrm{TV}} = 1.
\]
\end{theorem}
\begin{remark}
\label{remboundeddegree}
If $(G_n)$ is a sequence with $|V_n| \to\infty$ and \mbox{$\sup_{n} \ol
{\Delta}(G_n) < \infty$}, then Assumption \ref{assumpgraphs} is
equivalent to the decay of $g(x,y;G_n)$ in $d(x,y)$ uniformly in $n$.
\end{remark}

Many families satisfy Assumptions \ref{assumpgraphs} and \ref
{assumpfarawayharnack}, for example, $\Z_n^d$ for $d \geq3$, random
$d$-regular graphs whp, also for $d \geq3$, and the hypercube $\Z
_2^n$. We will discuss these and other examples in the next section.

The problem that we consider is closely related to determining the
mixing time of the \textit{lamplighter walk}, which we now introduce;
recall that $\CX(G) = \{ f \dvtx V \to\{0,1\}\}$ is the set of
markings of $V$ by $\{0,1\}$. If $G = (V,E)$ is a finite graph, the
wreath product $G^\diamond= \Z_2 \wr G$ is the graph $(V^\diamond,
E^\diamond)$ whose vertices are pairs $(f,x)$ where $f \in\CX(G)$
and $x \in V$. There is an edge between $(f,x)$ and $(g,y)$ if and only
if $\{x,y\} \in E$ and $f(z) = g(z)$ for $z \notin\{x,y\}$.
$G^\diamond$ is also referred to as the \textit{lamplighter graph} over
$G$ since it can be constructed by placing ``lamps'' at the vertices of
$G$; the first coordinate $f$ of a configuration~$(f,x)$ indicates the
state of the lamps and the second gives the location of the lamplighter.

%
\begin{figure}

\includegraphics{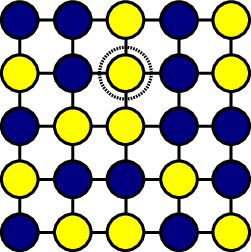}

\caption{A typical configuration of the lamplighter
over a $5 \times5$ planar grid. The colors indicate the state of the
lamps and the dashed circle gives the position of the lamplighter.}
\label{fig2}
\end{figure}

The \textit{lamplighter walk} $X^\diamond$ on $G$ is the random walk on
$G^\diamond$. Its transition kernel $p(\cdot,\cdot;G^\diamond)$ can
be constructed from $p(\cdot,\cdot;G)$ using the following procedure:
given $(f,x) \in V^\diamond$:
\begin{longlist}[(2)]
\item[(1)] sample $y \in V$ adjacent to $x$ using $p(x,\cdot;G)$,
\item[(2)] randomize the values of $f(x)$, $f(y)$ using independent fair
coin flips,
\item[(3)] move the lamplighter from $x$ to $y$.
\end{longlist}
See Figure \ref{fig2} for an example of a typical lamplighter configuration.
That both $f(x)$ and $f(y)$ are randomized rather than just $f(y)$ is
necessary for reversibility. It is obvious that the stationary
distribution of $X^\diamond$ is $\nu(\cdot;G) \times\pi(G)$. For
the graphs we consider, the mixing time of $X^\diamond$ is dominated
by the mixing time of its first coordinate as it is comparable to
$T_\cov(G)$ which in turn is much larger than $T_\mix(G)$, the mixing
time of the second coordinate of $X^\diamond$. This will\vspace*{2pt} allow us to
deduce $T_\mix(G^\diamond) = (\frac{1}{2}+o(1)) T_\cov(G_n)$ for
graphs satisfying Assumptions\vspace*{1pt} \ref{assumpgraphs}, \ref
{assumpfarawayharnack} from Theorem \ref{thmthreshold}.

Random walk on a sequence of graphs $(G_n)$ is said to have a (total
variation) cutoff with threshold $(a_n)$ if
\[
\lim_{n \to\infty} \frac{T_\mix(\varepsilon;G_n)}{a_n} = 1
\qquad\mbox{for all } \varepsilon\in(0,1).
\]
It is believed that many graphs have a cutoff, but establishing this is
often quite difficult since it requires a delicate analysis of the
behavior of the underlying walk. The term was first coined by Aldous
and Diaconis in \cite{AD86} where they prove cutoff for the
top-in-at-random shuffling process. Other early examples include random
transpositions on the symmetric group \cite{DS81}, the riffle shuffle
and random walk on the hypercube \cite{A83}. By making a small
modification to the proof of Theorem \ref{thmthreshold}, we are able
to establish cutoff for the lamplighter walk on base graphs satisfying
Assumptions \ref{assumpgraphs} and \ref{assumpfarawayharnack}.

Before we state these results, we will first summarize previous work
related to this problem. The mixing time of $G^\diamond$ was first
studied by H\"aggstr\"om and Jonasson in \cite{HJ97} in the case $G_n
= K_n$ and $G_n = \Z_n$. Their work implies a~cutoff with threshold
$\frac{1}{2} T_\cov(K_n)$ in the former case and that there is no
cutoff in the latter. The connection between $T_\mix(G^\diamond)$ and
$T_\cov(G)$ is explored further in \cite{PR04}, in addition to
developing\vspace*{1pt} the relationship between the relaxation time of $G^\diamond
$ and $T_\hit(G)$, and $\E[2^{|\CL(\alpha;G)|}]$ and $T_\mix
^U(G^\diamond)$. The results of \cite{PR04} include a proof of cutoff
when $G_n = \Z_n^2$ with threshold $T_\cov(\Z_n^2)$ and a~general
bound that
%
%
\begin{equation}
\label{eqnprbound}
\bigl[ \tfrac{1}{2} + o(1) \bigr] T_\cov(G_n) \leq T_\mix
(G_n^\diamond) \leq[1+o(1)] T_\cov(G_n),
\end{equation}
whenever $(G_n)$ is a sequence of vertex transitive graphs with $T_\hit
(G_n) = o(T_\cov(G_n))$. It is not possible to improve upon (\ref
{eqnprbound}) without further hypotheses since the lower and upper
bounds are achieved by $K_n$ and $\Z_n^2$, respectively.

The bound (\ref{eqnprbound}) applies to $\Z_n^d$ when $d \geq3$
since $T_\hit(\Z_n^d) \sim c_d n^d$ and $T_\cov(\Z_n^d) = c_d' n^d
(\log n)$ (see Proposition 10.13, Exercise 11.4 of \cite{LPW08}). This
leads \cite{PR04} to the question of whether there is a threshold for
$T_\mix((\Z_n^d)^\diamond)$ and, if so, if it is at $\frac{1}{2}
T_\cov(\Z_n^d)$, $T_\cov(\Z_n^d)$ or somewhere in between. By
a~slight extension of our methods, we are able to show that the threshold
is at $\frac{1}{2} T_\cov(\Z_n^d)$ when $d \geq3$, and that the
same holds whenever $(G_n)$ satisfies Assumptions \ref{assumpgraphs}
and \ref{assumpfarawayharnack}.\vspace*{-2pt}
\begin{theorem}
\label{thmlamplighter}
If $(G_n)$ satisfies Assumptions \ref{assumpgraphs} and \ref
{assumpfarawayharnack}, then $T_\mix(\varepsilon;\break G_n^\diamond)$ has a
cutoff with threshold $\frac{1}{2} T_\cov(G_n)$.\vspace*{-2pt}
\end{theorem}

In order to prove Theorems \ref{thmthreshold} and \ref
{thmlamplighter}, we need to develop a delicate understanding of both
the process of coverage and the correlation structure of~$\CL(\alpha
;G_n)$. The proof yields the following theorem, which gives a precise
estimate of the decay of correlation in $\CL(\alpha;G_n)$ under the
additional hypothesis of vertex transitivity.\vspace*{-2pt}
\begin{theorem}
\label{thmcorrelationdecay}
Suppose $(G_n)$ is a sequence of vertex transitive graphs satisfying
Assumption \ref{assumpgraphs}. If $(x_n^i)$ for $1 \leq i \leq\ell$
is a family of sequences with $x_n^i \in V_n$ and $|x_n^i - x_n^j| \geq
r$ for every $n$ and $i \neq j$, then
%
%
\begin{eqnarray}
\label{eqncorrelationdecay}
(1-\delta_{r,\ell}) |V_n|^{-\ell\alpha- \delta_{r,\ell}} &\leq&
\p[ x_n^i \in\CL(\alpha; G_n) \mbox{ for all } i]\nonumber\\[-8pt]\\[-8pt]
&\leq&(1+\delta_{r,\ell}) |V_n|^{-\ell\alpha+\delta_{r,\ell}},
\nonumber
\end{eqnarray}
where $\delta_{r,\ell} \to0$ as $r \to\infty$ while $\ell$ is
fixed. If $\ol{\Delta}(G_n) \to\infty$, we take $r = 1$ and $\delta
_{1,\ell} = o(1)$ as $n \to\infty$.\vspace*{-2pt}
\end{theorem}

\subsection*{Outline} The remainder of the article is structured as
follows. We show in Section \ref{secexamples} that the hypotheses
of\vadjust{\goodbreak}
Theorems \ref{thmthreshold} and \ref{thmlamplighter} hold for
a number of natural examples. In Section \ref{secprelim}, we collect
several general estimates that will be used throughout the rest of the
article; Proposition \ref{proptvbound} is in particular of critical importance.

Next, in Section \ref{seccovhit} we will develop precise asymptotic
estimates for the cover and hitting times of graphs $(G_n)$ satisfying
Assumption \ref{assumpgraphs}. The key idea is that the process by
which $X$ hits a point $x$ can be understood by studying the excursions
of $X$ from $\partial B(x,r)$ through $\partial B(x,R)$, $r < R$ and
then subsequently run for time $\beta T_\mix^U(G)$, some $\beta> 0$,
in order to remix. Due to the remixing, these excursions exhibit
behavior which is close to that of i.i.d. random walk excursions
initialized from stationarity. This has three important consequences.
First, our transience assumptions imply that the number $N_H(x)$ of
excursions up until the time $\tau(x)$ that $x$ is hit is
stochastically dominated from below by a geometric random variable with
small parameter $p$ provided $R > r$ are both large. Thus, $N_H(x)$ is
typically very large. Second and consequently, the empirical average of
the amount of time separating the beginning of successive excursions up
to time $\tau(x)$ is very concentrated around its mean $T_{r,R}(x)$.
Third, with $p_j(x)$ the probability that the $j$th excursion $E_j$
hits $x$ by time $\alpha T_\mix(G)$ after exiting $B(x,R)$, $\alpha
\leq\beta$, conditional on both the entrance point of $E_j$ and
$E_{j+1}$ to $B(x,r)$, we have that $\frac{1}{k}\sum_{j=1}^k p_j(x)$
is also well concentrated around its mean $\ol{p}_{r,R}(x)$. Combining
everything, this allows us to deduce the following asymptotic formula
for the hitting time of $x$:
\[
\E[ \tau(x)] = \bigl(1+o(1)\bigr) T_{r,R}(x) \E[ N_H(x)] = \frac{(1+o(1))
T_{r,R}(x)}{\ol{p}_{r,R}(x)}.
\]

For simplicity, we will now restrict our attention to graph families
which are vertex transitive. This implies that $T_{r,R} = T_{r,R}(x)$
and $\ol{p}_{r,R} = \ol{p}_{r,R}(x)$ do not depend on $x$.
Consequently, by the Matthews method upper and lower bounds (\cite
{M88}; see also Theorem 11.2 and Proposition 11.4 of \cite{LPW08}) we
infer that
%
%
\begin{equation}
\label{eqncovest}
T_\cov(G_n) = \bigl(1+o(1)\bigr) {\frac{T_{r,R}}{\ol{p}_{r,R}} \log}|V_n|.
\end{equation}

We will now explain how we use these estimates to prove Theorems \ref
{thmthreshold} and~\ref{thmlamplighter} in Section \ref
{sectotalvariation}. By Proposition \ref{proptvbound}, to give an
upper bound on the total variation distance of the i.i.d. marking of
the range of random walk run for time $\frac{1}{2} T_\cov(G_n)$ from
the uniform marking on $V_n$, it suffices to control the exponential
moment of the set of points in $V_n$ which are not visited by two
independent random walks, each run for time $\frac{1}{2} T_\cov
(G_n)$. Equation (\ref{eqncovest}) implies that the number $N_\cov
(x;\alpha)$ of excursions that have occurred by time $\alpha T_\cov
(G_n)$ satisfies
\[
N_\cov(x;\alpha) = {\bigl(\alpha+o(1)\bigr) \log}|V_n| / \ol{p}_{r,R}.\vadjust{\goodbreak}
\]
This in turn implies the tail decay
\[
\p[ \tau(x) \geq\alpha T_\cov(G_n)] = |V_n|^{-\alpha+ o(1)}.
\]
For points $x,y$ which are far apart, it is unlikely that a single
random walk excursion passes through both $B(x,R)$ and $B(y,R)$. That
is, the process of hitting well-separated points exhibits mean-field
behavior, which in turn allows us to give an efficient estimate of the
relevant exponential moment. There are many technical challenges
involved in getting all of these estimates to fit together correctly.

Decomposing the process of hitting into excursions between concentric
spheres is not new, and is used to great effect, for example, in \cite
{DPRMAN03,DPRZTHICK01,DPRZCOV04,DPRZLATE06}. Our implementation of this
idea is new since explicit representations of hitting probabilities and
Green's functions in addition to the approximate rotational invariance
available in the special case of $\Z_n^d$ are not available in the
generality we consider.

We prove Theorem \ref{thmcorrelationdecay} in Section \ref
{seccorrdecay}. This result, which may be of independent interest, is
important in Section \ref{sectotalvariation} since it allows us to
deduce that points in $\CL(\frac{1}{2}; G_n)$ are typically ``spread
apart.'' The article ends with a list of related open questions.

\section{Examples}
\label{secexamples}

\subsection*{$\Z_n^d$, $d \geq3$} Although the simplest, this is the
motivating example for this work. It is well known (see Section 1.5 of
\cite{LAW91}) that there exists a constant $c_d > 0$ so that $g(x,y;\Z
_n^d) \leq c_d |x-y|^{2-d}$, which implies uniform local transience.
Assumption~\ref{assumpfarawayharnack}\hyperlink{assumpharnack}{(2)} is
also satisfied since it is also a basic result that random walk on $\Z
_n^d$ satisfies a Harnack inequality (see \cite{LAW91}, Section 1.4).

\subsection*{Super-critical percolation cluster} Suppose that $\eta
_e$ is a collection of i.i.d. random variables indexed by the edges
$e=(x,y)$ of $\Z^d$, $d \geq3$, taking values in $\{0,1\}$ such that
$\p[\eta_e = 1] = p \in[0,1]$. An edge $e$ is called open if $\eta
_e = 1$. Let $\CC(x)$ denote the subset of $\Z^d$ consisting of those
elements $y$ that can be connected to $x$ by a path consisting only of
open edges. Let $C_\infty$ denote the event that there exists an
infinite open cluster and let $p_c = \inf\{ p > 0 \dvtx\break \p[ C_\infty]
>0\}$. Suppose $p > p _c$. Then it is known that there exists a unique
infinite open cluster $\CC_\infty$ almost surely. Fix $x \in\CC
_\infty$ and consider the graph $G_n = B(x,n) \cap\CC_\infty$. It
follows from the works of Delmotte~\cite{D99}, Deuschel and Pisztora
\cite{DP96}, Pisztora \cite{P96} and Benjamini and Mossel \cite{BM03}
that the heat kernel for continuous time random walk (CTRW) on
$G_n$ has Gaussian tails whp when $n$ is large enough; see the
discussion after the statement of Theorem A of \cite{B04}.
Consequently,  Green's function of the CTRW on $(G_n)$ has the same
quantitative behavior as for $(\Z_n^d)$. This implies the same is true
for the lazy random walk, which in turn yields uniform local transience
for $(G_n)$ whp when $n$ is sufficiently large. Therefore there exists
$n_0 = n_0(\omega)$ such that $(G_n\dvtx n \geq n_0(\omega))$ almost
surely satisfies Assumption \ref{assumpgraphs}.\vadjust{\goodbreak} Furthermore, it is a
result of Barlow \cite{B04} that there exists $n_1 = n_1(\omega)$
such that random walk on $(G_n\dvtx n \geq n_1(\omega))$ almost surely
satisfies a Harnack inequality and hence Assumption \ref{assumpfarawayharnack}.

\subsection*{Bounded degree expanders}
Suppose that $(G_n)$ is an expander family with uniformly bounded
maximal degree such that $|V_n| \to\infty$. Then there exists $T_0 <
\infty$ such that $T_{\mathrm{rel}}(G_n) \leq T_0$ for every $n$ where
$T_{\mathrm{rel}}(G_n)$ is the~re\-laxation time of lazy random walk on $G_n$.
Equation (12.11) of \cite{LPW08} implies that
\[
p^t(x,y;G_n) \leq C \biggl( \frac{1}{|V_n|} + e^{-t/T_0}\biggr)
\]
and Theorem 12.3 of \cite{LPW08} gives $T_\mix^U(G_n) = O({\log}
|V_n|)$. By Remark \ref{remboundeddegree}, to check Assumption \ref
{assumpgraphs}, we need only show the uniform decay $g(x,y;G_n)$ in~$d(x,y)$. If $t < d(x,y)$, then it is obviously true that $p^t(x,y;G_n)
= 0$. Hence,
%
%
\begin{eqnarray}
\label{eqnexpandergreendecay} g(x,y;G_n) &\leq& C\Biggl(\frac{ O({\log}
|V_n|)}{|V_n|} + \sum_{t=d(x,y)}^{T_\mix^U(G_n)} e^{- t/
T_0}\Biggr)\nonumber\\[-8pt]\\[-8pt]
&\leq& C_1 e^{- d(x,y) / T_0} + o(1)\nonumber
\end{eqnarray}
as $n \to\infty$.
We will now argue that $(G_n)$ satisfies part \hyperlink
{assumpfaraway}{(1)} of
Assumption~\ref{assumpfarawayharnack}. Suppose that $\ol{\Delta}
\geq\max_{x \in V_n} \deg(x)$ for every $n$. We can obviously take
$R_n^\gamma= {\gamma\log}|V_n| / (2 \log\ol{\Delta})$, hence we
have $T_\mix^U(G_n) / R_n^\gamma= O(1)$ as $n \to\infty$. Combining
this with (\ref{eqnexpandergreendecay}) implies that $(G_n)$ satisfies
Assumption \ref{assumpfarawayharnack}.

\subsection*{Random regular graphs} Suppose that $d \geq3$ and let
$\CG_{n,d}$ denote the set of $d$-regular graphs on $n$ vertices. It
is well known \cite{BS89} that, whp as \mbox{$n \to\infty$}, an element
chosen uniformly from $\CG_{n,d}$ is an expander. Consequently, whp, a
sequence $(G_n)$ where each $G_n$ is chosen independently and uniformly
from~$\CG_{n,d}$, $d \geq3$, almost surely satisfies the hypotheses
of our theorems.

\subsection*{Hypercube} As in the case of super-critical percolation,
for $\Z_2^n$ it is easiest to prove bounds for the CTRW which, as we
remarked before, easily translate over to the corresponding lazy walk.
The transition kernel of the CTRW is
\[
p^t(x,y;\Z_2^n) = \frac{1}{2^n} (1+e^{-2t/n})^{n-|x-y|} (1-e^{-2t/n})^{|x-y|},
\]
where $|x-y|$ is the number of coordinates in which $x$ and $y$ differ.
The spectral gap is $1/n$ (see Example 12.15 of \cite{LPW08}) which
implies $\Omega(n) = T_\mix^U(\Z_2^n) = O(n^2)$ (see Theorem 12.3
of
\cite{LPW08}). Suppose that $A \subseteq\Z_2^n$ has diameter $s$ and
$d(x,A) = r$. If $y \in A$, we have
\[
p^t(x,y;\Z_2^n) \leq\frac{1}{2^n} (1+e^{-2t/n})^{n-r} (1-e^{-2t/n})^r.
\]
It is easy to see that
\[
p^t(x,y;\Z_2^n) \leq\cases{
\displaystyle \biggl(C_\varepsilon\frac{t}{n}\biggr)^r \exp\biggl(- \frac{t}{C_\varepsilon n}
(n-r)\biggr), &\quad if $t \leq\varepsilon
n$,\vspace*{2pt}\cr
\displaystyle e^{-\rho_\varepsilon n}, &\quad if $t > \varepsilon n$,}
\]
provided $\varepsilon> 0$ is sufficiently small. Consequently,
\[
g(x,A;\Z_2^n) \leq C n^{s-r}
\]
and therefore $\Z_2^n$ is uniformly locally transient. The other
hypotheses of Assumption \ref{assumpgraphs} are obviously satisfied.
As for Assumption \ref{assumpfarawayharnack}, we note that in this
case, we can take $R_n^\gamma= \gamma n/ (2\log_2 n)$. Thus, if $r >
0$ it is easy to see that if $\diam(A) \leq s$ and $d(x,A) \geq
R_n^\gamma$ we have that
\[
\sum_{y \in A} p^t(x,y;\Z_2^n) \leq n^s e^{-\rho_\varepsilon n},
\]
if $t > \varepsilon n$. On the other hand, if $t \leq\varepsilon n$,
then we have
\[
\sum_{y \in A} p^t(x,y;\Z_2^n) \leq n^{s} \biggl( \frac{ C_\varepsilon
t}{n} \biggr)^{\gamma n/ (2\log_2 n)} e{- t/(2 C_\varepsilon)}.
\]
Hence, it is not hard to see that $\Z_2^n$ satisfies Assumption \ref
{assumpfarawayharnack}.

\subsection*{Caley graph of $S_n$ generated by transpositions} Let
$G_n$ be the Caley graph of $S_n$ generated by transpositions. By work
of Diaconis and Shahshahani \cite{DS81}, $T_\mix(G_n) = \Theta(n
(\log n))$, which by Theorem 12.3 of \cite{LPW08} implies $T_\mix
^U(G_n) = O( n^2 (\log n)^2)$. We are now going to give a crude
estimate of~$p^t(\sigma,\allowbreak\tau;S_n)$. By applying an automorphism, we
may assume without loss of generality that $\sigma= \mathrm{id}$. Suppose
that $d(\mathrm{id},\tau) = r$ and that $\tau_1,\ldots,\tau_r$ are
transpositions such that \mbox{$\tau_r \cdots\tau_1 = \tau$}. Then $\tau
_1,\ldots,\tau_r$ move at most $2r$ of the $n$ elements of $\{
1,\ldots,n\}$, say, $k_1,\ldots,k_{2r}$. Suppose $k_1',\ldots
,k_{2r}'$ are distinct from $k_1,\ldots,k_{2r}$ and $\alpha\in S_n$
is such that $\alpha(k_i) = k_i'$ for $1 \leq i \leq r$. Then the
automorphism of~$G_n$ induced by conjugation by $\alpha$ satisfies
$\alpha\tau\alpha^{-1} \neq\tau$. Therefore, the size of the set
of elements $\tau'$ in $S_n$ such that there exists a~graph
automorphism $\varphi$ of~$G_n$ satisfying $\varphi(\tau) = \tau'$
and $\varphi(\mathrm{id}) = \mathrm{id}$ is at least ${n-2r\choose2r}
\geq2^{-2r} n^{2r} ((2r)!)^{-1}$ assuming $n \geq8r$. Therefore,
\[
p^t(e,\tau;G_n) \leq\frac{2^{2r} (2r)!}{n^{2r}} \quad\mbox{and}\quad
g(e,\tau;G_n) \leq C(2^{2r} (2r)!) (\log n)^2 n^{2-2r}.
\]
If $\diam(A) = s$, then trivially $|A| \leq n^{2s}$ from which it is
clear that $(G_n)$ is uniformly locally transient. The other parts of
Assumption \ref{assumpgraphs} are obviously satisfied by $G_n$. As for
Assumption \ref{assumpfarawayharnack}, a simple calculation shows that
we can take $R_n^\gamma\leq\gamma n/4 + O(1)$. Hence setting
$R_n^\gamma= \sqrt{n}$, a calculation analogous to the one above,
gives that Assumption \ref{assumpfarawayharnack} is satisfied.

\section{Preliminary estimates}
\label{secprelim}

The purpose of this section is to collect several general estimates
that will be useful for us throughout the rest of the article.\vadjust{\goodbreak}

\begin{lemma}
\label{lemtvl1bound}
If $\mu,\nu$ are measures with $\nu$ absolutely continuous with
respect to $\mu$ and
\[
\int\frac{d\nu}{d\mu} \,d\nu= 1 + \varepsilon,
\]
then
\[
\|\nu- \mu\|_{\mathrm{TV}} \leq\frac{\sqrt{\varepsilon}}{2}.
\]
\end{lemma}
\begin{pf}
This is a consequence of the Cauchy--Schwarz inequality:
\begin{eqnarray*}
\|\mu- \nu\|_{\mathrm{TV}}^2
&=& \biggl( \frac{1}{2}\int\biggl| \frac{d\nu}{d\mu} - 1 \biggr|
\,d\mu\biggr)^2
\leq\frac{1}{4}\int\biggl| \frac{d\nu}{d\mu} - 1 \biggr|^2 \,d\mu\\[-2pt]
&=& \frac{1}{4}\biggl(\int\frac{d\nu}{d\mu} \,d\nu- 1\biggr).
\end{eqnarray*}
\upqed\end{pf}

Let $\nu$ denote the uniform measure on $\CX(G) = \{f \dvtx V \to\{
0,1\}\}$.\vspace*{-3pt}
\begin{proposition}
\label{proptvbound}
Suppose that $\mu$ is a measure on $\CX(G)$ given by first sampling
$\CR\subseteq V$ according to a probability $\mu_0$ on $2^{V}$, then,
conditional on $\CR$ sampling $f \in\CX(G)$ by setting
\[
f(x) = \cases{
\xi(x), &\quad if $x \in\CR$,\cr
0, &\quad otherwise,}
\]
where $(\xi(x) \dvtx x \in V)$ is a collection of i.i.d. random variables
with $\p[ \xi(x) = 0] = \p[\xi(x) = 1] = \frac{1}{2}$. Then
\[
\int\frac{d\mu}{d\nu} \,d\mu= \iint2^{|\CR^c \cap\CS^c|}
\,d\mu_0(\CR) \,d\mu_0(\CS).\vspace*{-3pt}
\]
\end{proposition}
\begin{pf}
Letting $\mu(\cdot|\CS)$ be the conditional law of $\mu$ given $\CS
$ and $N = |V|$, we have
\begin{eqnarray*}
\int\frac{d\mu}{d\nu} \,d\mu
&=& 2^N \int\mu(\{ f\}) \,d\mu(f)\\[-2pt]
&=& 2^N \iint\biggl( \int\mu(\{f\} | \CS) \,d\mu_0(\CS)
\biggr) \,d\mu(f|\CR) \,d\mu_0(\CR).
\end{eqnarray*}
Suppose $f \in\CX(G)$ is such that $f|_{\CR^c} \equiv0$ for some
$\CR\subseteq V$. Note that
\[
\mu(\{f\} | \CS) = 2^{-|\CR\cap\CS| - |\CS\setminus\CR|}| \one
_{\{f|_{\CR\setminus\CS\equiv0}\}}.
\]
Hence, the above is equal to
\begin{eqnarray*}
&&2^N \iint\biggl( \int2^{-|\CR\cap\CS|-|\CS\setminus\CR|}
\one_{\{f|_{\CR\setminus\CS\equiv0}\}} \,d\mu_0(\CS) \biggr) \,d\mu
(f|\CR) \,d\mu_0(\CR)\\[-2pt]
&&\qquad= 2^N \iint2^{-|\CR\cap\CS| - |\CS\setminus\CR|} \biggl(
\int\one_{\{f|_{\CR\setminus\CS\equiv0}\}} \,d\mu(f|\CR) \biggr)
\,d\mu_0(\CR) \,d\mu_0(\CS)\\[-2pt]
&&\qquad= 2^N \iint2^{-|\CR\cap\CS| - |\CS\setminus\CR|} 2^{-|\CR
\setminus\CS|} \,d\mu_0(\CS) \,d\mu_0(\CR).
\end{eqnarray*}
Simplifying the expression in the exponent gives the result.
\end{pf}

Roughly speaking, the general strategy of our proof will be to show
that if $\CR, \CR'$ denote independent copies of the range of random
walk on $G_n$ run up to time $(\frac{1}{2} + \varepsilon) T_\cov(G_n)$
and $\CL= V \setminus\CR$, $\CL' = V \setminus\CR'$ then
%
%
\begin{equation}
\label{eqnexpmoment}
\E\exp(\zeta|\CL\cap\CL'|) = 1 + o(1) \qquad\mbox{as } n \to\infty
\end{equation}
for $\zeta> 0$.
This method cannot be applied directly, however, since this exponential
moment blows up even in the case of $\Z_n^3$. To see this, suppose
that~$X,X'$ are independent random walks on $\Z_n^3$ initialized at
stationarity. We divide the cover time $c_3 n^3 (\log n)$ into rounds
of\vspace*{1pt} length $n^2$. In the first round, with probability $1/4$ we know
that $X$ starts in $\bL_1 = \Z_n^2 \times\{n/8,\ldots,3n/8\}$. In
each successive round, $X$ has probability $\rho_0 > 0$ strictly
bounded from zero in $n$ of not leaving $\bL_2 = \Z_n^2 \times\{
1,\ldots,n/2\}$ and ending the round in $\bL_1$. Since there are $c_3
n (\log n)$ rounds, this means that $X$ does not leave $\bL_1$ with
probability at least
\[
\tfrac{1}{4} \rho_0^{c_3 n \log n} \geq c \exp(-\rho_1 n \log n).
\]
Since $X'$ satisfies the same estimate, we therefore have
\[
\E\exp(\zeta|\CL\cap\CL'|) \geq c \exp\biggl(\frac{\zeta}{2} n^3-
2\rho_1 n\log n\biggr) \to\infty\qquad\mbox{as } n \to\infty.
\]
The idea of the proof is to truncate the exponential moment in (\ref
{eqnexpmoment}) by conditioning the law of random walk run for time
$(\frac{1}{2}+\varepsilon)T_\cov(G_n)$ \textit{conditional on typical
behavior} so that
\[
\| \wt{\mu}_0 - \mu_0 \|_{\mathrm{TV}} = o(1) \qquad\mbox{as } n \to\infty.
\]
We do this in such a way that the uncovered set exhibits a great deal
of spatial independence in order to make the exponential moment easy to
estimate. To this end, we will condition on two different events. The
first is that points in $\CL(\frac{1}{2} + \delta;G_n)$ are well separated: for any $x \in V_n$ the number of points in $\CL
(\frac{1}{2}+\varepsilon;G_n)$ which are contained in a large ball
centered at $x$ is at most some constant $M$.
Given this event, we can partition $\CL(\frac{1}{2} + \varepsilon
;G_n)$ into disjoint subsets $E_1,\ldots, E_M$ such that $x,y \in
E_\ell$ distinct implies $d(x,y)$ is large. Observe
\[
\E\exp(\zeta|\CL\cap\CL' \cap E_\ell|) \leq\E\prod_{x \in
E_\ell} \Biggl(1 + e^\zeta\prod_{j=1}^{N'(x,T)} \bigl(1-q_j'(x)\bigr)
\Biggr),
\]
where $N'(x,T)$ is the number of excursions of $X'$ from $\partial
B(x,r)$ to $\partial B(x,R)$ by time $T$ and $q_j'(x)$ is the
probability the $j$th such excursion hits $x$ conditional on its
entrance and exit points. When $T$ is large, uniform local transience
implies that $N'(x,T)$ and $\prod_{j=1}^k q_j'(x)$ can be\vspace*{1pt} estimated by
their mean and, roughly speaking, this is the second event on which we
will condition. Finally, we get control of the entire exponential
moment by an application of H\"older's inequality.

We finish the section by recording a standard lemma that bounds the
rate of decay of the total variation and uniform distances to stationarity:
\begin{proposition}
\label{propmixingdecay}
For every $s,t \in\N$,
%
%
\begin{eqnarray}
\label{eqntvdecay}\qquad
\max_x \| p^{t+s}(x,\cdot) - \pi\|_{\mathrm{TV}} &\leq& 4 \max_{x,y} \|
p^t(x,\cdot) - \pi\|_{\mathrm{TV}} \| p^s(y,\cdot) - \pi\|_{\mathrm{TV}} ,\\
\label{eqnuniformdecay}
\max_{x,y} \biggl| \frac{p^{t+s}(x,y)}{\pi(y)} - 1\biggr| &\leq&
\max_{x,y} \frac{p^{s}(x,y)}{\pi(y)} \max_x \| p^{t}(x,\cdot) -
\pi\|_{\mathrm{TV}}.
\end{eqnarray}
\end{proposition}
\begin{pf}
The first part is a standard result; see, for example, Lemmas~4.11
and~4.12 of \cite{LPW08}. The second part is a consequence of the
semigroup property:
\begin{eqnarray*}
\frac{1}{\pi(z)} p^{t+s}(x,z)
&=& \frac{1}{\pi(z)} \sum_y p^t(x,y) p^s(y,z)\\
&=& \frac{1}{\pi(z)} \sum_y [p^t(x,y) - \pi(y) + \pi(y)]p^s(y,z)\\
&\leq& \biggl(\max_{y,z} \frac{p^s(y,z)}{\pi(z)}\biggr) \|
p^t(x,\cdot) - \pi\|_{\mathrm{TV}} + 1.
\end{eqnarray*}
\upqed\end{pf}

Note that (\ref{eqntvdecay}) and (\ref{eqnuniformdecay}) give
%
%
\begin{eqnarray}\qquad
\max_{x} \| p^t(x,\cdot) - \pi\|_{\mathrm{TV}} &\leq& c e^{-c\alpha} \qquad
\mbox{for } t \geq\alpha T_\mix(G),\\
\label{eqnuniformexponentialdecay}
\max_{x,y} \biggl| \frac{p^{t+s}(x,y)}{\pi(y)} - 1 \biggr| &\leq& c
e^{-c\alpha} \qquad\mbox{for } t \geq T_\mix^U(G) + \alpha T_\mix(G),
\end{eqnarray}
where $c > 0$ is a universal constant. We will often use (\ref
{eqnuniformexponentialdecay}) without reference, and, for simplicity
use that the same inequality holds when $T_\mix^U(G) + \alpha T_\mix
(G)$ is replaced by $\alpha T_\mix^U(G)$, perhaps adjusting $c > 0$.

\section{Hitting and cover times}
\label{seccovhit}

Throughout, we assume that we have a sequence of graphs $(G_n)$
satisfying Assumption \ref{assumpgraphs} with transience function~$\rho$.
We will often suppress the index $n$ and refer to an element
of $(G_n)$ as $G$ and similarly write $V,E$ for $V_n,E_n$,
respectively. The primary purpose of this section is to develop
asymptotic estimates of the maximal hitting and cover times of $(G_n)$.
Roughly, these will be given in terms of:
\begin{longlist}[(2)]
\item[(1)] the return time $T_{r,R}(x)$, $x \in V$, of $X$ to $B(x,r)$ after
passing through $B(x,R)$, $R > r$, large then allowed to remix, and
\item[(2)] the probability $\ol{p}_{r,R}(x)$ that upon entering $B(x,r)$,
$X$ subsequently hits $x$ before exiting $B(x,R)$.
\end{longlist}
The derivation of these formulas requires many technical steps, so we
will provide an overview of how everything fits together before delving
into the details.

Let $N(x,t)$ be the number of excursions made by $X$ from $\partial
B(x,r)$\break to~$\partial B(x,R)$, then subsequently allowed to remix by
running for some multiple of $T_\mix^U(G)$, by time $t$ and let
$p_j(x)$ be the probability that the $j$th excursion~$E_j$ hits $x$
conditional on the entrance points of $E_j$ and $E_{j+1}$ to~$B(x,r)$.
Since the $p_j(x)$ are independent, we can express the probability~$P(x,t)$ that~$x$ has not been hit by time $t$ by the formula
\[
P(x,t) = \E\prod_{j=1}^{N(x,t)} \bigl(1-p_j(x)\bigr).
\]
We will argue using uniform local transience that we can make $p_j(x)$
as small as we like by choosing $R > r$ large enough. Consequently, we have
\[
P(x,t) = \E\Biggl[ \exp\Biggl(-\bigl(1+O(\rho(r))\bigr)\sum_{j=1}^{N(x,t)}
p_j(x)\Biggr) \Biggr].
\]

Our first goal, accomplished in the next subsection, is to show that
the empirical mean $\frac{1}{k} \sum_{j=1}^{k} p_j(x)$ is
concentrated around its mean $\ol{p}_{r,R}(x)$. Next, in Section
\ref{subsecexcursionlength}, we will again use concentration to argue
that $N(x,t) \approx t/T_{r,R}(x)$. These two steps allow us to
conclude that $P(x,t)$ is approximately given by $\exp(-t \ol
{p}_{r,R}(x) / T_{r,R}(x))$. That is, $P(x,t)$ is approximately
exponential with parameter $\ol{p}_{r,R}(x) / T_{r,R}(x)$ so that the
expected hitting time of $x$ is approximately $T_{r,R}(x) / \ol
{p}_{r,R}(x)$. In the vertex transitive case, this immediately leads to
an estimate of ${(T_{r,R} / \ol{p}_{r,R}) \log}|V|$ for the cover time
via the Matthews method (\cite{M88}; see also Theorem 11.2 and
Proposition 11.4 of \cite{LPW08}). A similar but more complicated
formula also holds for graphs which are not vertex transitive and is
derived in the second half of Section \ref{subsechittingcovering}.

\subsection{Probability of success}
\label{subsecprobsuccess}

Fix $R > r$ and let $X$ be a lazy random walk on $G$. Suppose $A = \{
x_1,\ldots,x_\ell\} \subseteq V$ where $d(x_i,x_j) \geq2R$ for $i
\neq j$. Let $A(s) = \{ x \in V \dvtx d(x,A) \leq s\}$ where $d(x,A) =
\min
_{y \in A} d(x,y)$. Let $\partial A(s) = \{x \in V\dvtx d(x,A) = s\}$. The
purpose of this section is to prove that the empirical mean of the
conditional probability that successive excursions of $X$ from~$\partial A(r)$ through $\partial A(R)$ succeed in hitting $x \in A$
given their entrance points concentrates around its mean. We will need
to extend our excursions by multiples of the uniform mixing time
$T_\mix^U(G)$ so we have enough independence to get good
concentration.\vadjust{\goodbreak}

%
\begin{figure}

\includegraphics{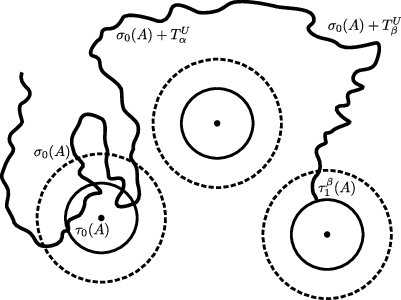}

\caption{The solid and dashed circles represent the
boundaries of $A(r)$ and $A(R)$, respectively, and the small points are
the elements of $A$. Note that $X$ may re-enter $A(r)$ during the
interval $[\sigma_k^\beta(A), \tau_{k+1}^\beta(A)]$.}
\label{fig3}
\end{figure}

To this end, we fix $\beta\geq0$, set $T_\beta^U = \beta T_\mix
^U(G)$, and define stopping times
%
%
\begin{eqnarray}
\label{covhiteqntau0}
\tau_0(A) &=& \min\{t \geq0\dvtx X(t) \in\partial A(r)\}, \\
\label{covhiteqnsig0}
\sigma_0(A) &=& \min\{ t \geq\tau_0(A)\dvtx X(t) \notin A(R)\}
\end{eqnarray}
and inductively set
%
%
\begin{eqnarray}
\label{covhiteqntauk}
\tau_k^\beta(A) &=& \min\{t \geq\sigma_{k-1}^{\beta}(A) + T_\beta
^U\dvtx X(t) \in\partial A(r)\}, \\
\label{covhiteqnsigk}
\sigma_k^\beta(A) &=& \min\{ t \geq\tau_k^\beta(A)\dvtx X(t) \notin
A(R)\} .
\end{eqnarray}
See Figure \ref{fig3} for an illustration of the stopping times described in
\mbox{(\ref{covhiteqntau0})--(\ref{covhiteqnsigk})}.
Fix $\alpha\in[0,\beta]$. Let $S_j^{\alpha,\beta}(x;A)$ be the
event that $X(t)$ hits $x$ in $[\tau_{j}^\beta(A),\break \sigma_{j}^\beta
(A) + T_\alpha^U]$,
\[
p_j^{\alpha,\beta}(x;A) = \p[S_j^{\alpha,\beta}(x;A)|X(\tau
_{j}^\beta(A)),X(\tau_{j+1}^\beta(A))]
\]
and
\[
a_j^{\alpha,\beta}(x;A) = \E\Biggl[\sum_{t = \tau_{j}^\beta
(A)}^{\sigma_{j}^\beta(A) + T_\alpha^U} \one_{\{X(t) = x\}} \Big|
X(\tau_{j}^\beta(A)) , X(\tau_{j+1}^\beta(A)) \Biggr].
\]
The reason that it is useful to consider $p_{r,R}^{\alpha,\beta
}(x;A)$ for $\beta> \alpha$ is that, as we will prove in Lemma \ref
{lemprobsuccessbound}, this allows us to show that the effect of
conditioning on the terminal point $X(\tau_{j+1}^\beta(A))$ of the
excursion is negligible when $\beta-\alpha$ is large enough. This in
turn allows us to use uniform local transience to get that
$p_{r,R}^{\alpha,\beta}(x;A)$ can be bounded in terms of the
transience function. Finally, we\vadjust{\goodbreak} let $\ol{p}{}^{\,\alpha,\beta
}_{\,r,R}(x;A) = \E_\pi p_0^{\alpha,\beta}(x;A)$ and $\ol{a}_{\,r,R}^{\,\alpha
,\beta}(x;A) = \E_\pi a_0^{\alpha,\beta}(x;A)$. For $\beta\geq
\alpha\geq1$ note that
\[
\ol{p}_{\,r,R}^{\,\alpha,\beta}(x;A) = \ol{p}_{r,R}^{1,1}(x;A) + O
\biggl( \frac{T_\beta^U}{|V|}\biggr)
\]
since a union bound implies that the probability $X$ hits $x$ in the
interval $[\sigma_0(A)+T_1^U, \sigma_0(A) + T_\beta^U]$ is $O(
T_\beta^U / |V|)$.

By Assumption \ref{assumpgraphs}, we have that $T_\beta^U \ol{\Delta
}{}^r(G) / |V| = o(1)$ as $n \to\infty$. Note that $\ol
{p}_{r,R}^{1,1}(x;A) \geq(2\ol{\Delta}(G))^{-r}$ since the
right-hand side\vspace*{1pt} bounds from below the probability that $X$ goes
directly from $\partial B(x,r)$ to $x$ in $r$ steps. Consequently,
%
%
\begin{equation}
\ol{p}_{\,r,R}^{\,\alpha,\beta}(x;A) = \bigl(1+o(1)\bigr) \ol{p}_{r,R}^{1,1}(x;A).
\end{equation}
From now on, we will write $\ol{p}_{r,R}(x;A)$ for $\ol
{p}_{r,R}^{1,1}(x;A)$. By the same argument, it is also true that $\ol
{a}_{\,r,R}^{\,\alpha,\beta}(x;A) = (1+o(1))\ol{a}_{r,R}^{1,1}(x;A)$ and
we will also wri\-te~$\ol{a}_{r,R}(x;A)$ for $\ol{a}_{r,R}^{1,1}(x;A)$.
\begin{lemma}
For each $\delta> 0$ there exists $\gamma_0 > 0$ such that for $\beta
-\alpha\geq\gamma_0$ and all $n$ large enough we have
%
%
\begin{eqnarray}
\label{lempjbound}
1-\delta&\leq&\frac{p_j^{\alpha,\beta}(x;A)}{\p[ S_j^{\alpha
,\beta}(x) | X(\tau_j^\beta(A))]} \leq1+\delta,\\
1-\delta&\leq&\frac{a_j^{\alpha,\beta}(x;A)}{ \E[ \sum_{t =
\tau_{j}^\beta(A)}^{\sigma_{j}^\beta(A) + T_\alpha^U} \one_{\{
X(t) = x\}} | X(\tau_{j}^\beta(A)) ]} \leq1+\delta.
\end{eqnarray}
In particular, $p_j^{\alpha,\beta}(x;A) \leq(1+\delta) \rho(r)$
and $a_j^\beta(x;A) \leq(1+\delta) \rho(0) \rho(r)$ where~$\rho$
is the transience function.
\end{lemma}
\begin{pf}
Note that
\begin{eqnarray*}
&& \p\bigl[ X\bigl(\sigma_j^\beta(A) + T_\alpha^U\bigr) = z | X(\tau_j^\beta(A))
= z_j, X(\tau_{j+1}^\beta(A)) = z_{j+1}\bigr]\\
&&\qquad= \frac{\p[ X(\sigma_j^\beta(A) + T_\alpha^U) = z, X(\tau
_j^\beta(A)) = z_j, X(\tau_{j+1}^\beta(A)) = z_{j+1}]}{\p[ X(\tau
_j^\beta(A)) = z_j, X(\tau_{j+1}^\beta(A)) = z_{j+1}]}\\
&&\qquad= \biggl(\frac{\p[ X(\tau_{j+1}^\beta(A)) = z_{j+1} | X(\sigma
_j^\beta(A) + T_\alpha^U) = z]}{\p[ X(\tau_{j+1}^\beta(A)) =
z_{j+1}| X(\tau_j^\beta(A)) = z_j]}\biggr) \\
&&\qquad\quad{}\times
\p\bigl[X\bigl(\sigma_j^\beta(A) + T_\alpha^U\bigr) = z| X(\tau
_j^\beta(A)) = z_j\bigr].
\end{eqnarray*}
Mixing considerations imply
\[
\p[ X(\tau_{j+1}^\beta(A)) = z_{j+1}| X(\tau_j^\beta(A)) = z_j] =
[ 1+ O(e^{-c\beta}) ] \p_\pi[ X(\tau_0(A)) = z_{j+1}]
\]
and
\begin{eqnarray*}
&&\p\bigl[ X(\tau_{j+1}^\beta(A)) = z_{j+1}| X\bigl(\sigma_j^\beta(A) +
T_\alpha^U\bigr) = z\bigr]\\
&&\qquad= \bigl[ 1+ O\bigl(e^{-c(\beta-\alpha)}\bigr) \bigr] \p_\pi[ X(\tau_0(A)) =
z_{j+1}].
\end{eqnarray*}
Consequently, if $\mu_j$ denotes the law of $X(\sigma_j^\beta(A) +
T_\alpha^U)$ conditional\break on $X(\tau_j^\beta(A))$ and $X(\tau
_{j+1}^\beta(A))$ and $\mu$ is the law of $X(\sigma_j^\beta(A) +
T_\alpha^U)$ but conditional only on $X(\tau_j^\beta(A))$, we have
$1-\delta\leq d\mu_j / d\mu\leq1+\delta$ when $\beta- \alpha$ is
large enough. Thus,
\begin{eqnarray*}
p_j^{\alpha,\beta}(x;A)
&=& \int\p\bigl[ S_j^{\alpha,\beta}(x)| X(\tau_j^\beta(A)), X\bigl(\sigma
_{j}^\beta(A) + T_\alpha^U\bigr) = z, X(\tau_{j+1}^\beta(A))\bigr] \,d\mu
_j(z)\\
&\leq& (1+\delta)\int\p\bigl[ S_j^{\alpha,\beta}(x)| X(\tau_j^\beta
(A)), X\bigl(\sigma_j^\beta(A) + T_\alpha^U\bigr) = z\bigr] \,d\mu(z)\\
&=& (1+\delta)\p[S_j^{\alpha,\beta}(x) | X(\tau_{j}^\beta(A))].
\end{eqnarray*}
The lower bound for $p_j^{\alpha,\beta}(x;A)$ and the bounds for
$a_j(x;A)$ are proved similarly.\vspace*{-3pt}
\end{pf}

In the next lemma,\vspace*{-1pt} we will prove the concentration of $p_j^{\alpha
,\beta}(x;A)$\break and~$a_j^{\alpha,\beta}(x;A)$. The proof consists of
three main steps. First,\vspace*{-1pt} the previous lemma allows us to replace
$p_j^{\alpha,\beta}(x;A)$ by $\p[S_j^{\alpha,\beta}(x;A) | X(\tau
_j(A))]$ and\vspace*{-1pt} likewise for~$a_j^{\alpha,\beta}(x;A)$. Roughly,\vspace*{-1pt} the
next step is to use a stochastic domination argument to show that we
can\vspace*{-1pt} replace $\p[S_j^{\alpha,\beta}(x;A) | X(\tau_j(A))]$ by i.i.d.
variables with law $\p[S_1^{\alpha,\beta}(x,A) | X(\tau_1(A))]$.
The result then follows by an application of Cram\'{e}r's theorem.
\begin{lemma}
\label{lemprobsuccessconcen}
Fix $r > 0$ and $\delta\in(0,1)$. There exists $\gamma_0 > 0$
depending only on $r,\delta$ such that for all $R \geq r$, $\beta-
\alpha\geq\gamma_0$ and $n$ large enough we have
%
%
\begin{eqnarray}
&& \p\Biggl[ \sum_{j=1}^k p_j^{\alpha,\beta}(x;A) \notin[1-\delta
,1+\delta] \ol{p}_{r,R}(x;A) k \Biggr]\nonumber\\[-10pt]\\[-10pt]
&&\qquad \leq4\exp\biggl(-\frac{C \delta^2 \ol
{p}_{r,R}(x;A)}{\rho(r)} k \biggr) \nonumber
\end{eqnarray}
and
%
%
\begin{eqnarray}
&& \p\Biggl[ \sum_{j=1}^k a_j^{\alpha,\beta}(x;A) \notin[1-\delta
,1+\delta] \ol{a}_{r,R}(x;A) k \Biggr]\nonumber\\[-10pt]\\[-10pt]
&&\qquad \leq4\exp\biggl(-\frac{C \delta^2 \ol
{a}_{r,R}(x;A)}{\rho(r)} k \biggr) ,\nonumber
\end{eqnarray}
where $C > 0$ is independent of $r,R,\delta$.\vadjust{\goodbreak}
\end{lemma}
\begin{pf}
Let $\mu$ be the measure on $\partial A(r)$ induced by the law of
$X(\tau_0(A))$ given that $X$ has a stationary initial distribution.
For each $\delta> 0$, let $\CM(\delta)$ be the set of measures $\nu
$ on $\partial A(r)$ which are uniformly mutually absolutely continuous
with respect to $\mu$ in the sense that
%
%
\begin{equation}
\label{eqnmeasurebound}
\max_{z \in\partial A(r)} \biggl|\frac{\nu(z)}{\mu(z)} - 1
\biggr| + \max_{z \in\partial A(r)} \biggl|\frac{\mu(z)}{\nu(z)} - 1
\biggr| \leq\delta.
\end{equation}
Let $\mu_y(z) = \p_y[ X(\tau^\gamma(A)) = z]$ where $\tau^\gamma
(A) = \min\{ t \geq T_\gamma^U\dvtx X(t) \in\partial A(r)\}$. Mixing
considerations imply that $\mu_y \in\CM( C e^{-C \gamma})$ for some
$C > 0$.
Fix $\delta> 0$, $\delta' < \delta/2$, and take $\beta- \alpha=
\gamma$ so large that $C e^{-C \gamma} \leq\delta'/2$. Let $\ol
{\mu}$, $\ul{\mu}$ be elements of $\CM(\delta'/2)$ such that $\p
[S_0^{\alpha,\beta}(x) | X(\tau_0(A)) = Z]$ where $Z \sim\ol{\mu
}, \ul{\mu}$ stochastically dominates from above and below,
respectively, all other choices in $\CM(\delta'/2)$. Assume that
$\gamma_0$ is chosen sufficiently large so that the previous lemma
applies for $\delta'/2$ when $n$ is sufficiently large.

Let $(U_j), (L_j)$ be i.i.d. sequences with laws $\p[S_0^{\alpha
,\beta}(x) | X(\tau_0(A)) = Z]$, $Z \sim\ol{\mu}, \ul{\mu}$,
respectively. With $\ol{U} = \E U_1$ and $\ol{L} = \E L_1$, obviously
\[
(1-\delta') \ol{p}_{r,R}(x;A) \leq\ol{L} \leq\ol{U} \leq
(1+\delta') \ol{p}_{r,R}(x;A).
\]
By construction, we can find a coupling of $U_j, L_j, p_j^{\alpha
,\beta}(x;A)$ so that
\[
L_j \leq p_j^{\alpha,\beta}(x;A) \leq U_j \qquad\mbox{almost surely for
all } j.
\]

Corollary 2.4.5 of \cite{DZ98} implies
\[
\E e^{\lambda U_1} \leq\frac{1}{2\rho(r)} \bigl( \ol{U} e^{2
\lambda\rho(r)} + 2\rho(r)-\ol{U} \bigr)
\]
hence Exercise 2.2.26 of \cite{DZ98} gives that the Fenchel--Legendre
transform $\Lambda^*$ of the law of $U_1$ satisfies
\[
\Lambda^*(u) \geq\wt{\Lambda}^*(u)\equiv\frac{u}{2\rho(r)} \log
\biggl(\frac{u}{\ol{U}}\biggr) + \biggl(1-\frac{u}{2\rho(r)}
\biggr) \log\biggl(\frac{1-u/(2\rho(r))}{1-\ol{U}/(2\rho(r))}\biggr).
\]
As
\[
\wt{\Lambda}^*(\ol{U}) = (\wt{\Lambda}^*)'(\ol{U}) = 0
\quad\mbox{and}\quad (\wt{\Lambda}^*)''(u) \geq\frac{1}{2\rho(r) u}
\]
we have
\[
\inf_{u \geq(1+\delta') \ol{U}} \Lambda^*(u) \geq\frac{1}{4\rho
(r) \ol{U}} (\delta')^2 \ol{U}{}^2 = \frac{(\delta')^2 \ol
{U}}{4\rho(r)},
\]
assuming $\delta' < 1$.
Consequently, Cram\'{e}r's theorem (Theorem 2.2.3, part (c), of \cite
{DZ98}) implies that
%
%
\begin{equation}
\label{eqnpconcen}
\p\Biggl[ \sum_{i=1}^k U_i \leq(1+\delta') \ol{U} k \Biggr] \geq
1-2\exp\biggl(- \frac{(\delta')^2 \ol{U} k}{4\rho(r)} \biggr).
\end{equation}
An analogous estimate also holds for $(L_i)$ with $\ol{U}$ replaced by
$\ol{L}$. The proof of concentration for the $a_j^{\alpha,\beta
}(x;A)$ is the same.
\end{pf}

\subsection{Excursion lengths}
\label{subsecexcursionlength}

We will make use of the same notation in this subsection as in the
previous. The main result is Lemma \ref{lemexcursionconcen}, which is
that the empirical average of successive excursion lengths\vspace*{1pt} $\tau
_{k+1}^\beta(A) - \tau_k^\beta(A)$ is exponentially concentrated
around its mean. The proof requires two auxiliary inputs. The first,
Lemma~\ref{lemradoncontrol}, is an estimate of the Radon--Nikodym
derivative of the law of random walk conditioned not to hit $A(r)$ with
respect to the stationary measure $\pi$. The second, Lemma \ref
{lemmeanexcursionlength}, gives that the mean length of an excursion
does not depend strongly on its starting point. Let $\tau(A) = \tau_0(A)$.
\begin{lemma}
\label{lemradoncontrol}
For $\alpha, s \geq0$ we have
\[
\p_y[ X(T_\alpha^U) = z | \CA] = \bigl[ 1+ O\bigl( e^{-c \alpha} + |A|
\rho(s,r)\bigr) + o(1) \bigr] \pi(z) \qquad\mbox{as } n \to\infty,
\]
where $\CA= \{ \tau(A) \geq T_\alpha^U, d(X(T_\alpha^U),A) \geq s\}$.
\end{lemma}
\begin{pf}
For $z \in V$ with $d(z,A) \geq s$, observe
%
%
\begin{eqnarray} \label{eqnradonbayes}
\p_y[X(T_\alpha^U) = z | \CA]
&=& \frac{\p_y[X(T_\alpha^U) = z, \tau(A) \geq T_\alpha^U]}{\p
_y[\CA]} \nonumber\\
&=& \frac{\p_y[\tau(A) \geq T_\alpha^U | X(T_\alpha^U) = z] \p
_y[X(T_\alpha^U) = z]}{\p_y[\CA]} \\
&=& \bigl(1+O(e^{-c\alpha})\bigr)\frac{\p_y[\tau(A) \geq T_\alpha^U |
X(T_\alpha^U) = z] \pi(z)}{\p_y[\CA]}.\nonumber
\end{eqnarray}
Fix $\alpha' < \alpha$. The idea of the proof is now to argue it is
unlikely for $\tau(A)$ to occur in the interval $[T_\alpha^U -
T_{\alpha'}^U, T_\alpha^U)$. This allows us to replace $T_\alpha^U$
above by $T_{\alpha}^U - T_{\alpha'}^U$ in (\ref{eqnradonbayes}).
This in turn allows us to use mixing considerations to deduce that
conditioning on $\{X(T_\alpha^U) = z\}$ has little effect on the
probability of $\{\tau(A) \geq T_\alpha^U - T_{\alpha'}^U\}$. We compute
\begin{eqnarray*}
&& \p_y[\tau(A) \geq T_\alpha^U | X(T_\alpha^U) = z]\\
&&\qquad= \p_y[\tau(A) \geq T_\alpha^U - T_{\alpha'}^U | X(T_\alpha^U) =
z] \\
&&\qquad\quad{}-\p_y[T_\alpha^U > \tau(A) \geq T_\alpha^U - T_{\alpha'}^U |
X(T_\alpha^U) = z].
\end{eqnarray*}
We have
\begin{eqnarray*}
&& \p_y[\tau(A) \geq T_\alpha^U - T_{\alpha'}^U | X(T_\alpha^U) = z]
\\
&&\qquad= 1-\frac{\p_y[ \tau(A) < T_\alpha^U - T_{\alpha'}^U, X(T_\alpha
^U) = z]}{\p_y[X(T_\alpha^U) = z]} \\
&&\qquad= 1-\frac{1+O(e^{-c\alpha})}{\pi(z)} \p_y[X(T_\alpha^U) = z |
\tau(A) < T_\alpha^U - T_{\alpha'}^U]\\
&&\qquad\quad\hspace*{15.3pt}{}\times \p_y[\tau(A) < T_\alpha^U -
T_{\alpha'}^U] \\
&&\qquad= \p_y[\tau(A) \geq T_\alpha^U - T_{\alpha'}^U] + O\bigl(e^{-c(\alpha
-\alpha')}\bigr).
\end{eqnarray*}

Note that
\begin{eqnarray*}
&& \p_y[T_\alpha^U > \tau(A) \geq T_\alpha^U - T_{\alpha'}^U |
X(T_\alpha^U) = z]\\
&&\qquad= \frac{1+ O(e^{-c\alpha})}{\pi(y)\pi(z)} \p_y[T_\alpha^U > \tau
(A) \geq T_\alpha^U - T_{\alpha'}^U , X(T_\alpha^U) = z] \pi(y).
\end{eqnarray*}
By reversing time, we see that this is equal to
\begin{eqnarray*}
&&\frac{1+ O(e^{-c\alpha})}{\pi(y)} \p_z[ \tau(A) \leq T_{\alpha
'}^U, d(X(t),A) > r \mbox{ for all } T_{\alpha'}^U < t \leq T_\alpha
^U, X(T_\alpha^U) = y]\\
&&\qquad\leq \frac{1+O(e^{-c\alpha})}{\pi(y)}\p_z[X(T_\alpha^U) = y|\tau
(A) \leq T_{\alpha'}^U]\p_z[\tau(A) \leq T_{\alpha'}^U]\\
&&\qquad= \bigl(1+O\bigl(e^{-c(\alpha-\alpha')}\bigr)\bigr)\p_z[\tau(A) \leq
T_{\alpha'}^U].
\end{eqnarray*}
A union bound along with uniform local transience implies this is of
order $O(|A|\rho(s,r) + o(1))$.
With $\CA_1 = \{d(X(T_\alpha^U),A) \geq s\}$,
\begin{eqnarray*}
\p_y[\CA] &=& \p_y[\tau(A) \geq T_\alpha^U, \CA_1]\\
&=& \bigl(\p_y[\tau(A) \geq T_\alpha^U - T_{\alpha'}^U | \CA
_1] - \p_y[T_\alpha^U > \tau(A) \geq T_\alpha^U - T_{\alpha'}^U
| \CA_1] \bigr) \p_y[ \CA_1]\\
&=& \p_y[\tau(A) \geq T_\alpha^U - T_{\alpha'}^U] + O\bigl(
e^{-c(\alpha-\alpha')}+ |A|\rho(s,r) +o(1) \bigr),
\end{eqnarray*}
the last line coming from a similar analysis as before.
Consequently,
\[
\frac{\p_y[\tau(A) \geq T_\alpha^U | X(T_\alpha^U) = z]}{\p
_y[\CA]}
= 1 + O\bigl( e^{-c(\alpha-\alpha')} + |A|\rho(s,r) +o(1)\bigr).
\]
Taking $\alpha' = \alpha/2$ gives the lemma.
\end{pf}

Let $\tau_k(A) = \tau_k^0(A)$, $\sigma_k(A) = \sigma_k^0(A)$, and
$T_{r,R}(A) = \E_\pi[\tau_1(A) - \tau_0(A)]$. We will now show that
mean excursion length does not depend too strongly on the starting
point of $X$. The idea is to argue that $X$ will typically run for some
multiple of the mixing time before getting close to $A$ provided it is
initialized sufficiently far away from $A$, then invoke the previous
lemma to replace the induced law on $V$ by $\pi$.
\begin{lemma}[(Mean excursion length)]
\label{lemmeanexcursionlength}
For every $r,\delta> 0$ there exists $R_0 > r$ such that $R \geq R_0$ implies
\[
(1-\delta) T_{r,R}(A) \leq\min_{y \notin A(R)} \E_y \tau_0(A)
\leq\max_{y \notin A(R)} \E_y \tau_0(A) \leq(1+\delta) T_{r,R}(A)
\]
for all $n$ large enough.
\end{lemma}
\begin{pf}
We have that
\[
\E_\pi[ \tau_1(A) - \tau_0(A)]
= \E_\pi[ \sigma_0(A) - \tau_0(A)] + \E_\pi[ \tau_1(A) - \sigma_0(A)].
\]
Obviously,
\[
\E_\pi[\sigma_0(A) - \tau_0(A)]
\leq\max_{y \in A(r)} \E_y \sigma_0(A)
\leq c T_\mix^U(G)
\]
for some $c > 0$ since in each interval of length $T_\mix^U(G)$,
random walk started in $A(r)$ has probability uniformly bounded from
below of leaving $A(R)$ provided $n$ is large enough. It is also
obvious that
\[
\min_{y \notin A(R)} \E_y \tau_0(A) \leq\E_\pi[ \tau_1(A) -
\sigma_0(A)] \leq\max_{y \notin A(R)} \E_y \tau_0(A).
\]
The previous lemma implies
\[
(1- \delta) \E_\pi[\tau_0(A)] \leq\E_y[\tau_0(A) | \CA] \leq
T_\alpha^U + (1+\delta) \E_\pi[ \tau_0(A)]
\]
for all $y \notin A(R)$ provided we choose $R,\alpha,s,n$ large enough
to accommodate our choice of $\delta$. Hence,
\[
(1-\delta)\E_\pi[\tau_0(A) ] \leq\E_y[\tau_0(A)] \leq(1+\delta
)\E_\pi[\tau_0(A) ]
\]
as it is not difficult to see that $T_\mix^U(G) = o(T_{r,R}(A))$ as $n
\to\infty$. Therefore
\[
\max_{y_1,y_2 \notin A(R)} \frac{\E_{y_1} \tau_0(A)}{\E_{y_2} \tau
_0(A)} \leq1+\delta,
\]
which proves the lemma.
\end{pf}

We end with the main result of the subsection, the concentration of the
empirical average of excursion lengths. The proof is an adaptation of
\cite{DPRMAN03}, Lemma 24, to our setting and is based on Kac's moment
formula (\cite{FP99}, Equation 6) for the first hitting time of a strong
Markov process along with the approximate i.i.d. structure of excursion lengths.
\begin{lemma}[(Concentration of excursions)]
\label{lemexcursionconcen}
For each $\beta\geq0$ and $r, \delta> 0$ there exists $R_0 > r$ such that
%
%
\begin{eqnarray}
\label{eqnexcursionconcenlower}
\p_{y}[ \tau_k^\beta(A) \leq(1-\delta) T_{r,R}(A) k
] &\leq& e^{-C \delta^2 k},\\
\label{eqnexcursionconcenupper}
\p_{y}[ \tau_k^\beta(A) \geq(1+\delta) T_{r,R}(A) k
] &\leq& e^{-C \delta^2 k}
\end{eqnarray}
for all $R \geq R_0$, $y \in V$ and $n$ large enough.
\end{lemma}
\begin{pf}
First of all, it follows from Lemma \ref{lemmeanexcursionlength} that
\[
\max_{y} \E_y[\tau_0(A)] \leq CT_{r,R}(A)
\]
for some $C > 0$ provided $R,n$ are sufficiently large. Consequently,
Kac's moment formula (see \cite{FP99}, Equation 6) for the first
hitting time of a strong Markov process implies for any $j \in\N$ we
have that
%
%
\begin{equation}
\label{eqnkacmoment}
\max_y \E_y[ (\tau_0(A))^j ] \leq j! c^j T_{r,R}^j(A)\vadjust{\goodbreak}
\end{equation}
for some $c > 0$. This implies that there exists $\lambda_0 > 0$ so that
\[
\max_{y} \E_y \exp[ \lambda\tau_0(A)/T_{r,R}(A)] < \infty
\qquad\mbox{for all } \lambda\in(0,\lambda_0).
\]
Using $\E[\sigma_0(A) - \tau_0(A)] = o(T_{r,R}(A))$, a similar
argument implies that, by possibly decreasing $\lambda_0$,
\[
\max_{y} \E_y \exp[\lambda\sigma_0(A)/T_{r,R}(A)] < \infty
\qquad\mbox{for all } \lambda\in(0,\lambda_0).
\]
Combining the strong Markov property with $T_\beta^U = o(T_{r,R}(A))$ yields
\[
\max_{y} \E_y \exp[\lambda\tau_k^\beta(A) / T_{r,R}(A)] < \infty
\qquad\mbox{for all } \lambda\in(0,\lambda_0).
\]
Let $R_0$ be large enough so that the previous lemma implies
\[
(1-\delta/2) T_{r,R}(A) \leq\min_{y \notin A(R)} \E_y \tau_0(A)
\leq\max_{y \notin A(R)} \E_y \tau_0(A) \leq(1+\delta/2) T_{r,R}(A)
\]
for $R \geq R_0$ and $n$ large enough.
We compute
\begin{eqnarray*}
\max_{y \notin A(R)} \E_y e^{-\theta\tau_0(A)}
&\leq&1 - \theta\min_{y \notin A(R)} \E_y \tau_0(A) + \theta^2 \max
_{y \notin A(R)} \E_y \tau_0^2(A)\\
&\leq&1 - \theta(1-\delta/2) T_{r,R}(A) + \zeta\theta^2 \\
&\leq&\exp
\bigl(\zeta\theta^2 - \theta(1-\delta/2) T_{r,R}(A)\bigr),
\end{eqnarray*}
where $\zeta= c T_{r,R}^2(A)$ for some $c > 0$.
Since $\tau_0(A) \geq0$, Chebychev's inequality leads to (\ref
{eqnexcursionconcenlower}). Indeed,
\begin{eqnarray*}
&&\p_{y}[ \tau_k^\beta(A) \leq(1-\delta) T_{r,R}(A) k ]\\
&&\qquad\leq\exp\bigl(\theta(1-\delta) T_{r,R}(A) k\bigr) \E_{y} e^{- \theta\tau
_k^\beta(A)}\\
&&\qquad\leq \exp\bigl(\theta(1-\delta)T_{r,R}(A) k\bigr)\Bigl[ \max_{y \notin
A(R)} \E_y e^{-\theta\tau_0(A)} \Bigr]^k\\
&&\qquad\leq \exp\bigl(\theta(1-\delta)T_{r,R}(A) k\bigr) \exp\bigl(\zeta\theta^2 k -
\theta(1-\delta/2) T_{r,R}(A) k\bigr).
\end{eqnarray*}
Taking
\[
\theta= \frac{\delta T_{r,R}(A)}{ c_1\zeta}
\]
we get that
\begin{eqnarray*}
&&\p_{y}[ \tau_k^\beta(A) \leq(1-\delta) T_{r,R}(A) k ]\\
&&\qquad\leq\exp\bigl(\zeta\theta^2 k -\theta T_{r,R}(A) k \delta/2\bigr)\\
&&\qquad\leq \exp\bigl(\zeta\delta^2 T_{r,R}^2(A) k / (c_1^2\zeta^2)-\delta^2
T_{r,R}^2(A) k / (2 c_1\zeta)\bigr)\\
&&\qquad\leq\exp(-c \delta^2 k)
\end{eqnarray*}
provided we take $c_1$ sufficiently large.\vadjust{\goodbreak}

To prove (\ref{eqnexcursionconcenupper}), we need to bound
\begin{eqnarray*}
&&\p_{y}[ \tau_k^\beta(A) \geq(1+\delta) T_{r,R}(A) k
]\\
&&\quad\leq \exp\bigl(-\theta(1+\delta) T_{r,R}(A) k \bigr) \Bigl( e^{ \theta
T_\beta^U} \max_{y} \E_y e^{\theta\tau_0(A)} \max_{y \in A(r)} \E
_y e^{\theta[\sigma_0(A) - \tau_0(A)]} \Bigr)^k.
\end{eqnarray*}
We again take
\[
\theta= \frac{\delta T_{r,R}(A)}{c_1 \zeta}
\]
with $c_1$ to be fixed shortly, and note that
\begin{eqnarray*}
\max_{y } \E_y e^{\theta\tau_0(A)}
&\leq&\bigl(1+o(1)\bigr) \max_{y \notin A(R)} \E_y e^{\theta\tau_0(A)}\\
&\leq& \exp\bigl(\theta(1+\delta/2) T_{r,R}(A) + c_2 \zeta\theta^2 + o(1)\bigr).
\end{eqnarray*}
Since $\max_{y \in A(r)} \E_y [\sigma_0(A) - \tau_0(A)] =
o(T_{r,R}(A))$ as $n \to\infty$, Kac's formula yields
\[
\max_{y \in A(r)} \E_y e^{\theta[\sigma_0(A) - \tau_0(A)]} =
1+o(1) \qquad\mbox{as } n \to\infty.
\]
Since $T_\beta^U = o(T_{r,R}(A))$ as $n \to\infty$ as well, we have
\begin{eqnarray*}
&& \p_{y}[ \tau_k^\beta(A) \geq(1+\delta) T_{r,R}(A) k]\\
&&\qquad\leq \exp\bigl(-\theta(1+\delta) T_{r,R}(A) k + \theta(1+\delta/2)
T_{r,R}(A) k + c_2 \zeta\theta^2 k + o(1) k\bigr)\\
&&\qquad\leq \exp\bigl(-\theta\delta T_{r,R}(A) k/2 + c_2 \zeta\theta^2 k +
o(1) k\bigr).
\end{eqnarray*}
Taking $c_1 > 0$ large enough gives the result.
\end{pf}

\subsection{Hitting and covering}
\label{subsechittingcovering}

The purpose of this subsection is to estimate the maximal hitting time
(Lemma \ref{lemhittingestimate}) and cover time (Lemma \ref{lemcovertime}).
\begin{lemma}[(Hitting time estimate)]
\label{lemhittingestimate}
For every $\delta> 0$ there exists $r_0$ such that for each $r \geq
r_0$ there is an $R_0 > r$ so that if $R \geq R_0$ the following holds.
If $A_n = \{x_{n1},\ldots,x_{n\ell}\} \subseteq V_n$ with $d(x_{ni},
x_{nj}) \geq2R$ for $i \neq j$ and $y_n \in V_n$ is such that
$d(x_{ni}, y_n) \geq2R$ for all $n$, then
%
%
\begin{eqnarray}
1- \delta&\leq&\liminf_{n \to\infty} \frac{\E_{y_n} \tau
(x_{ni})}{ T_{r,R}(A_n)/\ol{p}_{r,R}(x_{ni};A)}\\
&\leq&\limsup_{n \to\infty} \frac{\E_{y_n} \tau(x_{ni})}{
T_{r,R}(A_n)/\ol{p}_{r,R}(x_{ni};A)} \leq1+\delta.
\end{eqnarray}
\end{lemma}

As the proof of the lemma is long, we pause momentarily to highlight
the main steps. The primary tools will be the results from the previous
subsections.\vadjust{\goodbreak} The first ingredient (though we leave this to the end of
the proof) is to argue that it is unlikely for $X$ to hit a point
$x_{nk} \in A_n$ in the ``remixing'' intervals $[\sigma_k^\beta(A) +
T_\alpha^U, \sigma_k^\beta(A) + T_\beta^U]$. Once\vspace*{-1pt} we have
established this, it suffices to estimate the expectation of the first
time $\wt{\tau}(x_{nk})$ that $X$ hits $x_{nk}$ in $\bigcup_k [\tau
_k^\beta(A), \sigma_k^{\beta}(A) + T_\alpha^U]$ in place of the
expectation of $\tau(x_{nk})$. In particular, this implies that the
probability that $x_{nk}$ is first hit by the $(j+1)$st excursion is
well approximated by
\[
\E_{y_n} \Biggl[p_{j+1}^{\alpha,\beta}(x_{nk};A_n) \prod
_{i=1}^j\bigl(1-p_i^{\alpha,\beta}(x_{nk};A_n)\bigr) \Biggr].
\]
We now apply the concentration of the empirical mean of the
$p_j^{\alpha,\beta}(x;A)$ proved in Lemma \ref{lemprobsuccessconcen}
in order to replace the product with $\exp(-(1+\break O(\rho(r)))j \ol
{p}_{r,R}(x_{nk}$; $A_n))$, where we recall that $\rho$ is the transience
function. We conclude that the mean number of excursions required to
hit $x_{nk}$ is approximately $1/\ol{p}_{r,R}(x_{nk};A)$. The result
now follows by invoking Lemma~\ref{lemexcursionconcen}.

\begin{pf*}{Proof of Lemma \ref{lemhittingestimate}}
We will omit the indices $n$ and $i$ and just write $x$ for $x_{ni}$,
$y$ for $y_n$ and $A$ for $A_n$.
Fix $r$ sufficiently large so that $\rho(r) < \delta^2/100$. Recall
that $S_k^{\alpha,\beta}(x;A)$ is the event that $X$ hits $x$ in
$[\tau_k^\beta(A), \sigma_k^\beta(A) + T_\alpha^U]$ where $\tau
_k^\beta(A), \sigma_k^\beta(A)$ are as in (\ref
{covhiteqntau0})--(\ref{covhiteqnsigk}). Let $N(x;A) = \min\{ k \geq
1\dvtx S_k^{\alpha,\beta}(x;A)\break \mbox{occurs}\}$ and let
\[
\wt{\tau}(x) = \min\{ t \geq0\dvtx X(t) = x \mbox{ and } t \in I\},
\]
where
\[
I_k = [\tau_k^\beta(A), \sigma_k^\beta(A) + T_\alpha^U] \quad\mbox{and}
\quad I = \bigcup_k I_k.
\]
Then
\[
\tau_{N(x;A)}^\beta(A) \leq\wt{\tau}(x)
\leq\tau_{N(x;A)+1}^\beta(A).
\]
Let
\[
W(M;\delta) = \bigcap_{j \geq M} B(j;\delta) \equiv\bigcap_{j \geq
M} \{ (1-\delta) T_{r,R}(A) j \leq\tau_j^\beta(A) \leq
(1+\delta)T_{r,R}(A) j \}.
\]
With $\|\wt{\tau}(x) \| = \max_z \E_z \wt{\tau}(x)$, note that
%
%
\begin{eqnarray}\label{eqntildetautrunc}\quad
\E_y \wt{\tau}(x) \one_{W^c(M;\delta)}
&\leq&\sum_{j \geq M} \E_y \wt{\tau}(x) \one_{B^c(j;\delta)}
\nonumber\\
&\leq& \sum_{j \geq M} \bigl[\E_y \tau_j^\beta(x) \one
_{B^c(j;\delta)} + \| \wt{\tau}(x) \| \p_y[ B^c(j;\delta)] \bigr]
\nonumber\\[-8pt]\\[-8pt]
&\leq& 2C_0 \sum_{j \geq M} [ j T_{r,R}(A) + \| \wt{\tau}(x)\|
] e^{-C \delta^2 j} \nonumber\\
&\leq& C_1 \| \wt{\tau}(x)\| \sum_{j \geq M} (1+j) e^{-C \delta^2 j}
\leq C_2 \| \wt{\tau}(x)\| \frac{e^{-C \delta^2 M}}{\delta^4}.\nonumber
\end{eqnarray}
To see the second step, we let
\[
\wt{\tau}_j(x) = \min\{t \geq\tau_j^\beta(x)\dvtx X(t) = x\}.
\]
Then we have that
\begin{eqnarray*}
\E_y \wt{\tau}(x) \one_{B^c(j;\delta)} &\leq&
\E_y \wt{\tau}_j(x) \one_{B^c(j;\delta)} =
\E_y\bigl[ \bigl(\tau_j^\beta(x) + \bigl(\wt{\tau}_j(x) - \tau_j^\beta
(x)\bigr)\bigr) \one_{B^c(j;\delta)} \bigr]\\
&\leq&\E_y \tau_j^\beta(x) \one_{B^c(j;\delta)} + \E_y\bigl[\bigl(\wt{\tau
}_j(x) - \tau_j^\beta(x)\bigr) | B^c(j;\delta)\bigr] \p_y[B^c(j;\delta)].
\end{eqnarray*}
By the strong Markov property, $\E_y[\wt{\tau}_j(x) - \tau_j^\beta
(x) | B^c(j;\delta)] \leq\| \wt{\tau}(x)\|$.
In the third step, we used that
\begin{eqnarray*}
\E_y \tau_j^\beta(A) \one_{B^c(j;\delta)}
&\leq&(\E_y [\tau_j^\beta(A)]^2)^{1/2} (\p_y[B^c(j;\delta)])^{1/2}\\
&\leq&
\frac{2T_{r,R}(A)}{\lambda} j \bigl( \E_y \exp\bigl(\lambda\tau
_j^\beta(A) / (j T_{r,R}(A))\bigr) \bigr)^{1/2} C e^{-C \delta^2 j},
\end{eqnarray*}
where $\lambda\in(0,\lambda_0)$, $\lambda_0$ as in the proof of
Lemma \ref{lemexcursionconcen}. We used in the fourth step that
$T_{r,R}(A) = O(\| \wt{\tau}(x)\|)$. Indeed, this is true since
uniform local transience implies that with uniformly positive
probability more than one excursion is required to hit $x$ and, by
Lemma \ref{lemmeanexcursionlength}, the mean length of the second
excursion is at least $\frac{1}{2} T_{r,R}(A)$. The final step in
(\ref{eqntildetautrunc}) comes from summing the geometric series.
Uniform local transience implies
%
%
\begin{equation}
\label{eqnwttaumax}
\bigl| \E_y \wt{\tau}(x) - \| \wt{\tau}(x) \| \bigr| \leq\delta
\E_y \wt{\tau}(x),
\end{equation}
when $R$ is large enough. Consequently, there exists $M > 0$ large
enough depending only on $\delta$ so that
\[
\E_y \wt{\tau}(x) \one_{W(M;\delta)} \leq\E_y \wt{\tau}(x)
\leq(1+\delta) \E_y \wt{\tau}(x) \one_{W(M;\delta)}.
\]
Now,
\begin{eqnarray*}
\E_y \tau_{N(x;A)+1}^\beta(A) \one_{W(M;\delta)}
&=& \E_y \biggl[ N(x;A) \biggl( \frac{\tau_{N(x;A)+1}^\beta
(A)}{N(x;A)}\biggr) \one_{W(M;\delta)} \biggr]\\
&\leq& (1+\delta) T_{r,R}(x) \E_y N(x;A) +
\E_y \tau_M^\beta(A)\\
&\leq&(1+\delta) T_{r,R}(A) \E_y N(x;A) + C M T_{r,R}(A).
\end{eqnarray*}
In order to derive the inequality, we used that if $N(x;A) \geq M$ then
by the definition of $W(M;\delta)$ we have $\tau_{N(x;A)+1}^\beta(A)
/ N(x;A) \leq(1+\delta)T_{r,R}(A)$ and, in case $N(x;A) < M$,
we\vadjust{\goodbreak}
clearly have that $\tau_{N(x;A)+1}^\beta(A) \leq\tau_M^\beta(A)$.
The final\vspace*{1pt} inequality is a consequence of Lemma \ref
{lemmeanexcursionlength}. Similarly, we also have
\[
\E_y \tau_{N(x;A)}(A) \one_{W(M;\delta)} \geq(1-\delta)
T_{r,R}(A) \E_y N(x;A).
\]
Therefore,
\[
(1\,{-}\,\delta) T_{r,R}(A) \E_y N(x;A)\!\leq\!\E_y
\wt{\tau}(x)\!\leq\!(1\,{+}\,2\delta) T_{r,R}(A) \E_y N(x;A)\,{+}\,C M T_{r,R}(A).
\]
By Lemma \ref{lemprobsuccessconcen},
\begin{eqnarray*}
&&\ol{p}_{r,R}(x;A) \biggl[\exp\bigl(-(1+\delta) \ol{p}_{r,R}(x;A) j\bigr)-
C\exp\biggl(-\frac{C \delta^2 \ol{p}_{r,R}(x;A)}{\rho(r)} j
\biggr) \biggr]\\[-2pt]
&&\qquad\leq \E_y p_{j+1}^{\alpha,\beta}(x;A) \exp\Biggl(-[1+O(\rho
(r))]\sum_{i=1}^j p_i^{\alpha,\beta}(x;A) \Biggr)\\[-2pt]
&&\qquad\leq \ol{p}_{r,R}(x;A) \biggl[\exp\bigl(-(1-\delta) \ol{p}_{r,R}(x;A)
j\bigr) +C \exp\biggl(-\frac{C \delta^2 \ol{p}_{r,R}(x;A)}{\rho(r)}
j\biggr)\biggr].
\end{eqnarray*}
Taking $r$ sufficiently large gives
\begin{eqnarray*}
&&\E_y N(x;A)\\[-2pt]
&&\qquad= \sum_{j=1}^\infty j \p[ N(x;A) = j]\\[-2pt]
&&\qquad\leq CM^2 \rho(r) + \sum_{j=M+1}^\infty j \bigl(1+o(1)\bigr)\bigl(\ol
{p}_{r,R}(x;A) \exp\bigl(-(1-\delta) \ol{p}_{r,R}(x;A) j\bigr) \bigr)\\[-2pt]
&&\qquad\leq 2C M^2 \rho(r) + \frac{1+\delta}{\ol{p}_{r,R}(x;A)}.
\end{eqnarray*}
Similarly,
\[
\E_y N(x;A) \geq\frac{1-\delta}{\ol{p}_{r,R}(x;A)}.
\]
Increasing $r$ if necessary so that $M^2 \rho(r) \leq\delta$ yields
%
%
\begin{equation}
\label{eqnnxaestimate} \frac{1-2\delta}{\ol{p}_{r,R}(x;A)} \leq\E
_y N(x;A) \leq\frac{1+2\delta}{\ol{p}_{r,R}(x;A)}.
\end{equation}
This proves that
\[
\E_y \wt{\tau}(x) = \bigl(1+o(1)\bigr) \frac{T_{r,R}(A)}{\ol{p}_{r,R}(x;A)}
\qquad\mbox{as } n \to\infty.
\]
Let $F_k$ be the event that $X$ hits $A(r)$ in $J_k = [\sigma_k^\beta
(A) + T_\alpha^U, \sigma_k^\beta(A) + T_\beta^U]$. With $F = \bigcup
_{k=1}^{N(x;A)+1} F_k$, we have
\[
\E_y \wt{\tau}(x) \one_{F^c} \leq\E_y \tau(x) \leq\E_y \wt
{\tau}(x),
\]
where we recall that $\tau(x)$ is the first time $X$ hits $x$.\vadjust{\goodbreak}

We now claim that
%
%
\begin{equation}
\label{eqntautruncest}\qquad
\E_y \wt{\tau}(x) \one_{F^c} = \biggl[ 1 + O\biggl( \frac{ T_\beta
^U |A| \ol{\Delta}{}^r(G) }{|V| \ol{p}{}^2_{r,R}(x;A)} + \ol
{p}_{r,R}^2(x;A) \biggr)^{1/2} \biggr]\E_y \wt{\tau}(x).
\end{equation}
Note that this will complete the proof of the lemma as $\ol
{p}_{r,R}(x;A) \geq C\ul{\Delta}^{-r}(G)$ so that, by Assumption \ref
{assumpgraphs}, the error term can be made as small as we like by
making $r,R$ large enough. Using the Kac moment formula (\cite{FP99},
Equation~6) in the second inequality, we trivially have
%
%
\begin{eqnarray}\label{eqnwttaufbound}
\E_y \wt{\tau}(x) \one_F
&\leq& \E_y \tau(x) \one_F + \| \wt{\tau}(x)\| \p[F]
\nonumber\\[-8pt]\\[-8pt]
&\leq& C_1 \| \tau(x)\| \sqrt{\p[F]} + \| \wt{\tau}(x)\| \p[F].\nonumber
\end{eqnarray}
In view\vspace*{1pt} of (\ref{eqnwttaumax}) we have $\| \tau(x) \|\!\leq\!\| \wt
{\tau}(x)\|\!\leq\!(1\,{+}\,\delta) \E_y \wt{\tau}(x)$. Thus, using $\p
[F] \leq\break\sqrt{\p[F]}$, we see that we can bound (\ref
{eqnwttaufbound}) from above by $C_2 \| \wt{\tau}(x)\| \sqrt{\p
[F]}$. Using exactly the same proof of (\ref{eqnnxaestimate}), we have that
%
%
\begin{equation}
\label{eqnnxaestimatem2} \E_y[ N^2(x;A)] \leq\frac{C_3}{\ol
{p}_{r,R}^2(x;A)}.
\end{equation}
Applying (\ref{eqnnxaestimatem2}) along with Markov's inequality in
the second step, we consequently have
\begin{eqnarray*}
\p_y[F]
&\leq& \p_y[F, N(x;A)+1 \leq1/(\ol{p}_{r,R}(x;A))^2] \\
&&{}+ \p[N(x;A)
+1\geq1/(\ol{p}_{r,R}(x;A))^{2}]\\
&\leq& \sum_{k=1}^{1/(\ol{p}_{r,R}(x;A))^{2}} \p_y[F_k] + O((\ol
{p}_{r,R}(x;A))^2).
\end{eqnarray*}
Since $|A(r)|!\leq\!|A| \ol{\Delta}{}^r(G)$, a union bound implies $\p
_y[F_k]\!=\!O(T_\beta^U |A| \ol{\Delta}{}^r(G) / |V|)$, which proves
(\ref{eqntautruncest}).
\end{pf*}

If $G$ were vertex transitive so that $\ol{p}_{r,R}(x)$ and
$T_{r,R}(x)$ did not depend on~$x$, then by the Matthews method (\cite
{M88}; see also Theorem 11.2 and Proposition 11.4 of \cite{LPW08}) it
is possible to deduce that $T_{\cov}(G)$ is asymptotically well
approximated by $T_{r,R} / {\ol{p}_{r,R} \log}|V|$. Our goal\vspace*{1pt} is to
prove something similar even if $G$ is not vertex transitive. The idea
of the proof will be to group vertices together based on their hitting
time $T_{r,R}(x) / \ol{p}_{r,R}(x)$. In particular, we will argue that
the amount of time it takes to cover a set $V_F \subseteq V$ of
vertices each of whose hitting time is close $T_F$ is approximately
${T_F \log}|V_F|$. The cover time of $G$ is then well approximated by
${\max_{F} T_F \log}|V_F|$ where~$F$ ranges over subsets of vertices
with approximately constant hitting time.

The first step in implementing this strategy is to show that if we want
to estimate $T_\cov(G)$ to a multiple of $\varepsilon T_\cov(G)$,
$\varepsilon> 0$ fixed, we only need to consider a finite
number,\vadjust{\goodbreak}
depending only on $\varepsilon$, of groups of vertices. This will be
accomplished by relating $\ol{p}_{r,R}(x) / T_{r,R}(x)$ to $\pi(x)$
and then invoking Assumption \ref{assumpgraphs}.

We will now specialize to the case $A = \{x\}$; for simplicity of
notation we will omit $A$. Let
\[
O_{r,R}(x) = \frac{\ol{a}_{r,R}(x)}{T_{r,R}(x)} .
\]

\begin{lemma}
\label{lemoccestimate}
For every $\delta> 0$, there exists $r_0$ such that if $r \geq r_0$
there is $R_0 > r$ such that $R \geq R_0$ implies
\[
(1-\delta)\pi(x) \leq O_{r,R}(x) \leq(1+\delta)\pi(x)
\]
for all $n$ large enough.
\end{lemma}
\begin{pf}
Let $N(x,T) = \min\{ k\dvtx\tau_k^\beta(x) \geq T\}$, $J_k$ as in the
previous lemma, $J = \bigcup_k J_k$ and $\CG(x) = \sigma( X(\tau
_j^\beta(x))\dvtx j \geq1)$. Then
\[
\sum_{j=1}^{N(x,T)} a_j^{\alpha,\beta}(x) \leq\E\Biggl[ \sum
_{t=1}^T \one_{\{X(t) = x\}} \one_{\{t \notin J\}} \Big| \CG(x)
\Biggr] \leq\sum_{j=1}^{N(x,T)+1} a_j^{\alpha,\beta}(x).
\]
Lemmas \ref{lemprobsuccessconcen} and \ref{lemexcursionconcen} give that
\[
(1-\delta) T_{r,R}(x) \leq\frac{N(x,T)}{T} \leq(1+\delta)
T_{r,R}(x)
\]
and
\[
(1-\delta) \ol{a}_{r,R}(x) \leq\frac{\sum_{j=1}^k a_j^{\alpha
,\beta}(x)}{k} (1+\delta) \ol{a}_{r,R}(x)
\]
with high probability as $T \to\infty$, for all $r,R,k,n,\beta
-\alpha$ large enough.
Consequently, using that $(a_j^{\alpha,\beta}(x)\dvtx j \geq1)$ is
uniformly bounded, it is not hard to see that
\[
(1-\delta)\frac{\ol{a}_{r,R}(x)}{T_{r,R}(x)} \leq\frac{1}{T} \sum
_{j=1}^{N(x,T)} a_j^{\alpha,\beta}(x) \leq(1+\delta)\frac{\ol
{a}_{r,R}(x)}{T_{r,R}(x)}
\]
with high probability as $T \to\infty$, for all $r,R,n,\beta-\alpha
$ large enough. The middle term converges to $\pi(x)$ as $T \to\infty
$ since
\[
\lim_{T \to\infty} \frac{1}{T} \E\sum_{t=1}^T \one_{\{ X(t) \in
A(r) \}} \one_{\{t \in J\}} = 0.
\]
\upqed\end{pf}

Uniform local transience implies that there exists constants $c,C > 0$
so that $c \ol{a}_{r,R}(x) \leq\ol{p}_{r,R}(x) \leq C \ol
{a}_{r,R}(x)$; combining this with the previous lemma yields
\[
\frac{c \deg(x)}{|E|} \leq\frac{\ol{p}_{r,R}(x)}{T_{r,R}(x)} \leq
\frac{C \deg(x)}{|E|}.\vadjust{\goodbreak}
\]
Let $\varepsilon> 0$ and let
\[
H_{n,k}^\varepsilon= \biggl\{ x \in V_n\dvtx\frac{\ul{\Delta}(G_n) k
\varepsilon}{|E_n|} < \frac{\ol{p}_{r,R}(x)}{T_{r,R}(x)} \leq\frac
{\ul{\Delta}(G_n)(k+1)\varepsilon}{|E_n|} \biggr\}
\]
be a partition of $V_n$ into at most $\Delta_0 \varepsilon^{-1} $
subsets, where $\Delta_0$ is the constant from Assumption \ref
{assumpgraphs}. By passing to a subsequence, we may assume without loss
of generality that
\[
d_k^\varepsilon= \lim_{n \to\infty} d_{n,k}^\varepsilon\equiv\lim_{n
\to\infty} \frac{{\log}|H_{n,k}^\varepsilon|}{{\log}|V_n|}
\]
exists for every $k$. Note that $d_k^\varepsilon\in[0,1]$ for those $k$
so that $|H_{n,k}^\varepsilon| \neq0$ for all~$n$ large enough and,
since the partition is finite, necessarily there exists $k$ so that
$d_k^\varepsilon= 1$. In particular, there exists $k$ so that
$d_k^\varepsilon\neq0$. Let
%
%
\begin{equation}
C_{n,k}^\varepsilon= \frac{|E_n|}{ \ul{\Delta}(G_n) k \varepsilon}
{d_k^\varepsilon\log}|V_n| \quad\mbox{and}\quad
C_n^\varepsilon= \max_k C_{n,k}^\varepsilon.
\end{equation}
\begin{lemma}[(Cover time estimate)]
\label{lemcovertime}
For each $\delta> 0$, there exists $r_0,\varepsilon_0$ so that if $r
\geq r_0$ there is $R_0 > r$ such that $R \geq R_0$ and $\varepsilon\in
(0,\varepsilon_0)$ implies
%
%
\begin{equation}
\label{eqnsubsetcoverbound}
1-\delta\leq\liminf_{n \to\infty} \frac{T_\cov(H_{n,k}^\varepsilon
)}{C_{n,k}^\varepsilon} \leq\limsup_{n \to\infty} \frac{T_\cov
(H_{n,k}^\varepsilon)}{C_{n,k}^\varepsilon} \leq1+\delta
\end{equation}
for all $k$ with $d_k^{\varepsilon} > 0$.
Furthermore,
%
%
\begin{equation}
\label{eqncovertimebound}
1-\delta\leq\liminf_{n \to\infty} \frac{T_\cov
(G_n)}{C_n^\varepsilon} \leq\limsup_{n \to\infty} \frac{T_\cov
(G_n)}{C_n^\varepsilon} \leq1+\delta.
\end{equation}
\end{lemma}
\begin{pf}
Suppose $k$ is such that $d_k^\varepsilon> 0$. Then
$|H_{n,k}^\varepsilon|
\to\infty$ as $n \to\infty$. Let $r,R,n > 0$ be sufficiently large
so that Lemma \ref{lemhittingestimate} applies with our choice of
$\delta$. By Assumption \ref{assumpgraphs}\hyperlink
{assumpgraphsball}{(1)} we have that ${\log}|B(x,r)| = o({\log}|V_n|)$.
Consequently, for all $n$ large enough there exists an $R$-net
$E_{n,k}^\varepsilon$ of $H_{n,k}^\varepsilon$ such that
\[
{\log}|E_{n,k}^\varepsilon| = {\log}|H_{n,k}^\varepsilon| + o(1)
\qquad\mbox{as }
n \to\infty.
\]
The upper and lower bounds from the Matthews method (\cite{M88}; see
also Theorem~11.2 and Proposition 11.4 of \cite{LPW08}) combined with
the definition of~$C_{n,k}^\varepsilon$ imply~(\ref
{eqnsubsetcoverbound}). Theorem 2 of \cite{A91} implies that
\[
\lim_{n \to\infty} \frac{\tau_{\cov}(H_{n,k}^\varepsilon)}{\E\tau
_{\cov}(H_{n,k}^\varepsilon)} = 1.
\]
As $\tau_\cov(G_n) = \max_k \tau_\cov(H_{n,k}^\varepsilon)$ and the
maximum is over a finite set, it follows that $\tau_\cov(G_n) =
(1+o(1))\max_k T_\cov(H_{n,k}^\varepsilon)$. Taking expectations of
both sides gives (\ref{eqncovertimebound}).
\end{pf}

\section{Correlation decay}
\label{seccorrdecay}

The purpose of this section is to prove Theorem~\ref
{thmcorrelationdecay}. Exactly the same proof will also yield Lemma
\ref{lemcorrdecayupperbound}, a technical result which will be useful
in the next section, which is stated after the proof. Note that vertex
transitivity implies $\ol{p}_{r,R}(\cdot)$ and $T_{r,R}(\cdot)$ do
not depend on their arguments.
\begin{pf*}{Proof of Theorem \ref{thmcorrelationdecay}}
First, assume that we are in the case of bounded maximal degree. Let
$A$ be as in the previous section and let $\delta> 0$ be arbitrary.
Fix $r$ so that $\rho(r) \leq\delta^3/100 C \ell$ where $\ell=
|A|$ and $\ol{p}_{r,R}(x;A) \leq\delta^3$ for all $x \in A$. Let
$R_0 > r$ and $\beta-\alpha$ be sufficiently large so that Lemmas
\ref{lemprobsuccessconcen} and~\ref{lemexcursionconcen} apply with
our choice of $\delta,r$. Finally, let $N(x_i;A) = \min\{ k\dvtx
S_k^{\alpha,\beta}(x_i;A)\break \mbox{occurs}\}$ and $\CG(A) = \sigma
(p_j^{\alpha,\beta}(x;A)\dvtx x \in A, j \geq1)$. Since $d(x_i,x_j)
\geq2R$, the probability that $X$ neither hits $x$ nor $x'$ in the
interval $[\tau_j^\beta(x;A), \sigma_j^\beta(x;A) + T_\alpha^U]$ is
%
%
\begin{equation}
\label{eqnpjAunionbound}
1 - [1+ O(\rho(R))][ p_j^{\alpha,\beta}(x;A) + p_j^{\alpha,\beta}(x';A)].
\end{equation}
Indeed, the reason for this is that the conditional probability $X$
hits $B(x',R)$ in the same excursion that it hits $x$ given that it
hits the latter first is~$O(\rho(R))$ and the probability that $X$
hits $x$ before $B(x',R)$ is trivially bounded by $p_j^{\alpha,\beta}(x;A)$.
This holds more generally for any subset of $A$, hence
%
%
\begin{eqnarray} \label{eqnnohitseveral}
&&\E\bigl[ \p[ N(x_1;A) > k_1,\ldots, N(x_\ell;A) > k_\ell| \CG
(A)] \bigr] \nonumber\\
&&\qquad= \E\prod_{i=1}^\ell\exp\Biggl( -[1 + O(\rho(R))] \sum
_{j=1}^{k_i} p_j^{\alpha,\beta}(x_i;A) \Biggr) \nonumber\\[-8pt]\\[-8pt]
&&\qquad= \exp\Biggl( -[1+O(\delta)] \sum_{i=1}^\ell\ol{p}_{r,R}(x_i;A)
k_i \Biggr) \nonumber\\
&&\qquad\quad{}+ \sum_{i=1}^\ell O\bigl(\exp\bigl(- \ol{p}_{r,R}(x_i;A) k_i /
\delta\bigr)\bigr),\nonumber
\end{eqnarray}
where the last equality followed from our choice of $r$ and Lemma \ref
{lemprobsuccessconcen}.
Let $J_k = [\sigma_k^\beta(A) + T_\alpha^U, \sigma_k^\beta(A) +
T_\beta^U]$, as the in the previous section.
Combining this with Lemma \ref{lemexcursionconcen} and that the
probability $X$ hits $A(r)$ in $J_k$ is at most $O(T_\beta^U|A| \ol
{\Delta}{}^r(G) / |V|) = o(\ol{p}_{r,R}(x;A))$ for any $x \in A$, we have
\begin{eqnarray*}
&& \p[ \tau(x_1) \geq k T_{r,R}(A) / \ol{p}_{r,R}(x_1;A) ,\ldots,
\tau(x_\ell) \geq k T_{r,R}(A) / \ol{p}_{r,R}(x_n;A)]\\
&&\qquad= \bigl(1+o(1)\bigr) \exp\bigl( -[1+O(\delta)] \ell k \bigr) + O\bigl(\exp\bigl( -
C\delta^2 k /\rho(r)\bigr)\bigr)\\
&&\qquad= \bigl(1+o(1)\bigr) \exp\bigl( -[1+O(\delta)] \ell k \bigr).
\end{eqnarray*}
By vertex transitivity,
\[
T_\hit(G) = \bigl(1+o(1)\bigr) \frac{T_{r,R}(x_i;A)}{\ol{p}_{r,R}(x_i;A)}.
\]
By Lemma \ref{lemcovertime}, we know that the cover time is
asymptotically ${T_\hit(G) \log}|V|$. Inserting this into (\ref
{eqnnohitseveral}) gives the result for bounded degree.

This proof works also for unbounded degree, but is not quite sufficient
for the statement of our theorem since we would like to allow for
points in~$A$ to be adjacent. There are two parts that break down.
First, in Section \ref{seccovhit} we proved the concentration of
$p_j^{\alpha,\beta}(x;A)$ when $x \in A$ and we also\vspace*{2pt} assumed that
$x,y \in A$ implies $d(x,y) \geq2R$. To allow for $x,y$ adjacent, we define
\[
p_j^{\alpha,\beta}(y;A) = \p[ S_j^{\alpha,\beta}(y;A) | X(\tau
_j^\beta(A)), X(\tau_{j+1}^\beta(A))]
\]
for $y \in A(r/2)$. It is not difficult to see that for such $y$,
$p_j^{\alpha,\beta}(y;A)$ exhibits nearly the same concentration
behavior as for $y \in A$. Second, the estima\-te~(\ref
{eqnpjAunionbound}) is no longer good enough since $\rho(1)$ does not
decay in $n$. However, it is not difficult to see that the same
probability satisfies the estimate
%
%
\begin{equation}
1 - [1+ O(\ol{\Delta}{}^{-1}(G)) ][ p_j^{\alpha,\beta}(x;A) +
p_j^{\alpha,\beta}(x';A)],
\end{equation}
which suffices since $\ol{\Delta}{}^{-1}(G_n) \to0$ as $n \to\infty
$. The rest of the proof is the same.
\end{pf*}

Vertex transitivity was used only to get that $T_{r,R}(x;A) / \ol
{p}_{r,R}(x;A) = (1+o(1)) T_\hit(G)$. The same proof works more
generally, but leads to more complicated formulas. However, it is not
difficult to see that the upper bound takes a very similar form. This
result will be especially useful in the next section to show that
points which have not been visited by $X$ after time~$\frac{1}{2}
T_\cov(G)$ are typically well separated. Precisely, our estimate is:
\begin{lemma}
\label{lemcorrdecayupperbound}
If $(x_n^i)$ for $1 \leq i \leq\ell$ is a family of sequences with
$x_n^i \in H_{n,k(i)}^\varepsilon$ and $|x_n^i - x_n^j| \geq r$ for every
$n$ and $i \neq j$,
%
%
\begin{equation}
\label{eqncorrelationdecayupperbound}
\p[ x_n^i \in\CL(\alpha; G_n) \mbox{ for all } i]
\leq(1+\delta_{r,\ell}) |V_n|^{-\ell d_k^\varepsilon\alpha+\delta
_{r,\ell}},
\end{equation}
where $\delta_{r,\ell} \to0$ as $r \to\infty$ while $\ell$ is
fixed. If $\ol{\Delta}(G_n) \to\infty$ then we take $r = 1$ and
$\delta_{1,\ell} = o(1)$ as $n \to\infty$.
\end{lemma}

\section{Total variation bounds}
\label{sectotalvariation}

We are now in a position to complete the proof of Theorems \ref
{thmthreshold} and \ref{thmlamplighter}. We will prove the
lower bound first since it does not require us to specialize depending
on whether $(G_n)$ satisfies part~\hyperlink{assumpfaraway}{(1)} or
\hyperlink{assumpharnack}{(2)} of Assumption~\ref
{assumpfarawayharnack}. As we have
explained earlier, the upper bound will be proved by estimating the
exponential moment of the set of points not visited by two independent
random walks $X,X'$, each run for time $\frac{1}{2} T_\cov(G)$. We
will use Lemma \ref{lemcorrdecayupperbound} in the proof of Lemma \ref
{lemwellseparated} to argue that those\vadjust{\goodbreak} points~$\CL$ not visited by $X$
are typically far apart. This will be useful very useful because, as we
prove in Section \ref{subsecqconcen}, the hypothesis of Assumption~%
\ref{assumpfarawayharnack} allows us to establish concentration for
the empirical average of the conditional probability $q_j(x)$ that
excursions between $\partial B(x,r)$ to $\partial B(x,R)$ given both
the entry and exit points, where $R > r$ are very large.

\subsection{Lower bound}

We will now prove the lower bound for Theorems \ref{thmthreshold}
and~\ref{thmlamplighter}. This is actually just a slight extension of
Theorem 4.1 of \cite{PR04}, but we include it for the reader's
convenience. Recall from the \hyperref[sec1]{Introduction} that $\mu
(\cdot;\alpha,G)$
is the probability measure on $\CX(G) = \{f \dvtx V \to\{0,1\}\}$
given by first sampling $\CR(\alpha;G)$ then setting
\[
f(x) = \cases{
\xi(x), &\quad if $x \in\CR(\alpha;G)$,\cr
0, &\quad otherwise,}
\]
where $(\xi(x)\dvtx x \in V)$ is a collection of i.i.d. variables such
that $\p[\xi(x) = 0] = \p[\xi(x) = 1] = \frac{1}{2}$ and $\nu
(\cdot;G)$ is the uniform measure on $\CX(G)$.
\begin{lemma}[(Lower bound)]
\label{lemlowerbound}
For every $\delta> 0$,
\[
\lim_{n \to\infty} \bigl\| \mu\bigl(\cdot;\tfrac{1}{2}-\delta, G_n\bigr) - \nu
(\cdot;G_n)\bigr\|_{\mathrm{TV}} = 1.
\]
\end{lemma}
\begin{pf}
For $A \subseteq V$ and $m > 0$, let $\tau_{\cov}(A;m)$ be the first
time all but $m$ of the vertices of $A$ have been visited by $X$. For
each $k$ such that $d_k^\varepsilon> 0$, we will show that
%
%
\begin{equation}
\label{eqncoverbound} \lim_{n \to\infty} \p[ \tau_\cov(
H_{n,k}^\varepsilon; |H_{n,k}^\varepsilon|^{\alpha}) < (1-\alpha-\delta)
C_{n,k}^\varepsilon] = 0
\end{equation}
for each $\delta> 0$ and $\varepsilon\in(0,\varepsilon_0(\delta))$.
If not, then for some such $k,\delta,\alpha$ we have
\[
\limsup_{n \to\infty} \p[ A_{n,k}(\alpha,\delta)] > 0,
\]
where
\[
A_{n,k}(\alpha,\delta) = \{\tau_\cov( H_{n,k}^\varepsilon;
|H_{n,k}^\varepsilon|^{\alpha}) < (1-\alpha-\delta) C_{n,k}^\varepsilon
\}.
\]
It follows from the Matthews method upper bound (\cite{M88}; see also
Theorem~11.2 of \cite{LPW08}) that
\begin{eqnarray*}
&&\E[ \tau_\cov(H_{n,k}^\varepsilon) - \tau_{\cov}(H_{n,k}^\varepsilon;
|H_{n,k}^\varepsilon|^\alpha) | A_{n,k}(\alpha,\delta)]\\
&&\qquad\leq \alpha\bigl(1+O(\varepsilon)\bigr) C_{n,k}^\varepsilon\leq\alpha(1+\delta
/4) C_{n,k}^\varepsilon,
\end{eqnarray*}
where we take $\varepsilon$ so small that the $O(\varepsilon)$ term is at
most $\delta/4$.
Markov's inequality now implies
\[
\p[ \tau_{\cov}(H_{n,k}^\varepsilon) < (1 - \delta/2)
C_{n,k}^\varepsilon| A_{n,k}(\alpha,\delta)] > 0.
\]
This is a contradiction as Theorem 2 of \cite{A91} implies $\tau
_{\cov}(H_{n,k}^\varepsilon)/ C_{n,k}^\varepsilon\to1$ in
probability.

For each $n$ let $k_0(n)$ be an index that achieves the maximum in
$\max_k C_{n,k}^\varepsilon$. Now, (\ref{eqncoverbound}) implies that
whp\vspace*{1pt} at time $\frac{1}{2}(1-3\delta) T_\cov(G_n) = \frac
{1}{2}(1-3\delta+\break O(\varepsilon)) C_{n,k_0(n)}^\varepsilon$ the size of\vadjust{\goodbreak}
the subset of $H_{n,k_0(n)}^\varepsilon$ not visited by $X$ is at least
$|H_{n,k_0(n)}^\varepsilon|^{(1+2\delta+ O(\varepsilon))/2}$ but less than
$|H_{n,k_0(n)}^\varepsilon|^{(1+4\delta+ O(\varepsilon))/2}$. Thus, the
number of zeros in a marking of $H_{n,k_0(n)}^\varepsilon$ sampled from
$\mu(\cdot;\frac{1}{2}(1-3\delta),G_n)$ is whp at least
\[
\tfrac{1}{2}\bigl| H_{n,k_0(n)}^\varepsilon\bigr| +
\bigl(1+o(1)\bigr)\bigl|H_{n,k_0(n)}^\varepsilon
\bigr|^{(1+2\delta+ O(\varepsilon))/2} \qquad\mbox{as } n \to\infty.
\]
This proves the lemma since the probability of having deviations of
this magnitude from the mean tends to zero in a uniform marking.
\end{pf}

\subsection{Concentration of $q_j$}
\label{subsecqconcen}

\begin{figure}
\begin{tabular}{cc}

\includegraphics{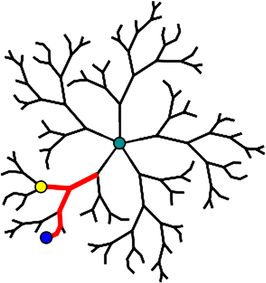}
 & \includegraphics{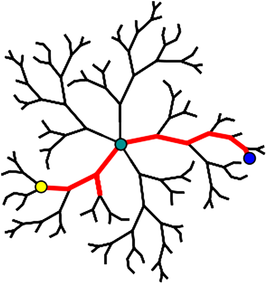} \\
(a) & (b)
\end{tabular}
\caption{The concentration\vspace*{-2pt} behavior of the $q_j(x)$ is very different
from the $p_j^{\alpha,\beta}(x)$ since it is not in general true that
$q_j(x) \leq C\rho(r)$ while it is true that $p_j^{\alpha,\beta}(x)
\leq C \rho(r)$. For example,\vspace*{1pt} in a graph which is locally tree like as
depicted above, it can be that $q_j(x) = 1$ for some combinations of
entrance and exit points.
\textup{(a)} Entrance and exit points of an excursion from $B(x,4)$ to $B(x,6)$,
respectively, conditional on which random walk has a~low probability of
hitting $x$.
\textup{(b)} Entrance and exit points of an excursion from $B(x,4)$ to~$B(x,6)$,
respectively, conditional on which random walk is forced to hit~$x$.}
\label{figqconcen}
\end{figure}

Let\vspace*{-1pt} $\sigma_j(x) = \sigma_j^{0,0}(x)$ and define $\tau_j(x)$
likewise where $\sigma_j^{\alpha,\beta},\tau_j^{\alpha,\beta}$
are as in (\ref{covhiteqntau0})--(\ref{covhiteqnsigk}). Let\vspace*{1pt} $S_j(x)$
be the event that $X$ hits~$x$ in the interval $[\tau_j(x),\sigma
_j(x)]$ and set $q_j(x) = \p[S_j(x) | X(\tau_j(x)), X(\sigma
_j(x))]$. The purpose of this subsection is to study the concentration
behavior of~$q_j(x)$, which will in turn depend on whether we assume
part \hyperlink{assumpfaraway}{(1)} or \hyperlink{assumpharnack}{(2)}
of Assumption
\ref{assumpfarawayharnack}; note that $q_j(x)$ differs from
$p_j^{\alpha,\beta}(x)$ from Section \ref{seccovhit}. Indeed, the
excursions on which we condition are different since we do not allow
the random walk to run for a multiple for $T_\mix^U(G)$ after exiting
$\partial B(x,R)$ and we\vspace*{1pt} condition on the entrance and exit points of
the current excursion rather than the entrance points of the current
and successive excursion. While both of these changes may seem
cosmetic, they affect the concentration behavior, since while
$p_j^{\alpha,\beta}(x)$ satisfies (\ref{lempjbound}), in locally
tree-like graphs it can be that $q_j(x) = 1$ with positive probability;
see Figure \ref{figqconcen} for an illustration of this behavior.

%

We shall first suppose that $(G_n)$ satisfies Assumption \ref
{assumpfarawayharnack}\hyperlink{assumpfaraway}{(1)}.
Let $\varepsilon> 0$ be arbitrary, $R_n^\gamma$ be as in Assumption \ref
{assumpfarawayharnack}, $\gamma> 0$ to be determined later, and let~%
$A$ be a set of points in $V_n$ such that if $x,y$ are distinct in $A$
then $d(x,y) \geq4 R_n^\gamma$. Fix $R > r > 0$ and let $\tau
_{k+1}(A) = \min\{ t \geq\sigma_k(x)\dvtx X(t) \in\partial A(r)\}$.
Fix $\beta> 0$ and define indices $i(j,x)$ inductively as follows. Set
\[
i(1,x) = \min\{ k \geq1\dvtx\tau_{k+1}(A) - \sigma_k(x) \geq T_\beta
^U \}
\]
and, for each $j \geq1$, let
\[
i(j+1,x) = \min\{ k \geq i(j,x) + 1 \dvtx
\tau_{k+1}(A) - \sigma_k(x) \geq T_\beta^U \}.
\]
When $x$ is clear from the context we will write $i(j)$ for $i(j,x)$.
\begin{lemma}
\label{lemiidreplacement}
For each $\delta> 0$ and $r > 0$ there exists $R_0 > r$ such that for
$R > R_0$ fixed there exists i.i.d. random variables $( I(j,x) \dvtx x
\in
A, j \geq1)$ which stochastically dominate from above $(i(j,x)\dvtx x
\in
A, j \geq1)$ and satisfy
\[
\p\bigl[ I\bigl((1-\delta)j,x\bigr) \geq j\bigr] \leq C\exp(-C \delta^2 j)
\]
for all $n$ large enough. Let $\CG(j,x) = \sigma(\{q_{i(k)}(x)\dvtx k
\neq j\} \cup\{ q_{i(k)}(y)\dvtx y \in A \setminus\{ x\}\})$.
There
exists i.i.d. random variables $(Q_j(x)\dvtx j \geq1)$ taking values in
$[0,2\rho(r)]$ such that
\[
1-O(e^{-c \beta}) \leq\frac{\E[ q_{i(j)} (x) | \CG(j,x)]}{Q_j(x)}
\leq1+O(e^{-c\beta})\vadjust{\goodbreak}
\]
and
\[
1-O(e^{-c \beta}) \leq\frac{\ol{p}_{r,R}(x)}{\E Q_j(x)} \leq1+
O(e^{-c \beta})
\]
for all $n$ large enough. Furthermore, the families $\{ (Q_j(x)\dvtx j
\geq1)\dvtx x \in A\}$ are independent.
\end{lemma}
\begin{pf}
Define stopping times
\begin{eqnarray*}
\sigma_{k0}(A) &=& \min\{ t \geq\sigma_k(x)\dvtx d(X(t),A) \geq
2R_n^\gamma\},\\
\tau_{k1}(A) &=& \min\{ t \geq\sigma_{k0}(x)\dvtx d(X(t),A) \leq
R_n^\gamma\}.
\end{eqnarray*}
For $j \geq1$, inductively set
\begin{eqnarray*}
\sigma_{k j }(A) &=& \min\{ t \geq\tau_{kj}(A)\dvtx d(X(t),A) \geq
2R_n^\gamma\},\\
\tau_{k(j+1)}(A) &=& \min\{ t \geq\sigma_{kj}(A)\dvtx d(X(t),A) \leq
R_n^\gamma\}.
\end{eqnarray*}
Note that $\sigma_{kj}(A) - \tau_{kj}(A) \geq R_n^\gamma$. Thus, for
$j_\beta= T_\beta^U / R_n^\gamma$ we have that $\tau_{kj_\beta}(A)
\geq\sigma_k(x) + T_\beta^U$. Let $F_k(x) = \{ X(t) \in A(r) \mbox{
for } t \in[\sigma_k(x), \sigma_k(x) + T_\beta^U]\}$. Let $x_{kj}$
be the element in $A$ such that $d(X(\tau_{kj}(A)), x_{kj}) \leq
R_n^\gamma$. Observe
\begin{eqnarray*}
&&\p_{X(\tau_{kj}(A))}\bigl[ X(t) \in A(r) \mbox{ for } t \in[ \tau
_{kj}(A), \sigma_{kj}(A)] | x_{kj}\bigr] \\
&&\qquad\leq C\max_{d(y,x_{kj}) = R_n^\gamma} g(y, B(x_{kj}, r);G_n).
\end{eqnarray*}
Uniform local transience also yields
\[
\p_{X(\sigma_k(x))}\bigl[ X(t) \in A(r) \mbox{ for } t \in[\sigma_k(x),
\tau_{k0}(A)]\bigr] \leq C \rho(R,r) \leq\delta/2,
\]
provided $R > r$ is large enough. A union bound thus gives
%
%
\begin{eqnarray}\label{eqnhitnetbound}
\p_{X(\tau_k(x))}[ F_k(x) ] &\leq& \max_{z} \max_{d(y,z) =
R_n^\gamma}g(y, B(z, r);G_n)\frac{T_\beta^U}{R_n^\gamma} + \delta
/2 \nonumber\\[-8pt]\\[-8pt]
&\leq& \delta/2 + o(1)\leq\delta\nonumber
\end{eqnarray}
as $n \to\infty$ by part \hyperlink{assumpfaraway}{(1)} of Assumption
\ref
{assumpfarawayharnack}. Note that if $x_1,\ldots,x_\ell\in A$ and
$j(1),\ldots,j(k)$ are such that $\tau_{j(k)}(x_k) \leq\tau
_{j(k+1)}(x_{k+1})$ then we have
\begin{eqnarray*}
&&\p\bigl[ F_{j(1)}(x_1),\ldots, F_{j(\ell)}(x_\ell)\bigr]\\
&&\qquad= \E\bigl[ \p_{X(\tau_{j(\ell))}(x_\ell)}\bigl[ F_{j(\ell)}(x_\ell)\bigr] \one
_{F_{j(1)}(x_1)} \cdots\one_{F_{j(\ell-1)}(x_{\ell-1})}\bigr]\\
&&\qquad\leq \delta\p\bigl[ F_{j(1)}(x_1),\ldots, F_{j(\ell-1)}(x_{\ell-1})\bigr]
\leq\cdots\leq\delta^{\ell}.
\end{eqnarray*}
This can of course be repeated with any subset of the above events
which implies the stochastic domination claim. It easily now follows
from Cram\'{e}r's theorem that
\[
\p\bigl[ I\bigl((1-\delta)k,x\bigr) \geq k \bigr] \leq2 \exp(-C \delta^2 k).
\]

For the second part of the lemma, we just need to get a bound on $\mu
_x(z) / \pi(z)$ where $\mu_x$ is the law of random walk started at
$x$ conditioned not to get within distance $r$ of $A$ by, say, time
$T_{\beta/2}^U$. This can be done in exactly the same way as in the
proof of Lemma \ref{lemradoncontrol}. Indeed, the term $|A| \rho
(s,r)$ in the statement of that lemma comes from a bound on the
probability that random walk at distance $s$ from $A$ hits $A$ in time
$T_\alpha^U$. In the situation of this lemma, the role of $s$ is
replaced by $R_n^\gamma$ and we can use the scheme developed above to
estimate the error contributed by this term by $O(\delta)$ provided
$n$ is sufficiently large.
\end{pf}

We now turn to the case that $(G_n)$ satisfies part \hyperlink
{assumpharnack}{(2)}
of Assumption \ref{assumpfarawayharnack}. This case
will turn out to be substantially easier, the reason being that the
Harnack inequality implies the quenched bound $q_j(x) \leq2C \rho
(r)$. We emphasize once more that this is not the case in locally
tree-like graphs.
\begin{lemma}
\label{lemprobsuccessbound}
If $(G_n)$ satisfies part \hyperlink{assumpharnack}{(2)} of Assumption
\ref
{assumpfarawayharnack}, then for each $r,\delta> 0$ there exists $R_0
> r$ such that $R \geq R_0$ implies
%
%
\begin{eqnarray}
\label{eqnqconcentrationharnack}
&&\p\Biggl[ \prod_{j=1}^{k} \bigl(1-q_j(x)\bigr) \geq\bigl(1-(1+\delta) \ol
{p}_{r,R}(x)\bigr)^{k(1+\delta)}\Biggr]\nonumber\\[-8pt]\\[-8pt]
&&\qquad \leq C\bigl[ \exp\bigl( - C \delta^2 \ol
{p}_{r,R}(x) k / \rho(r) \bigr) + \exp(- C \delta^2 k) \bigr]
\nonumber
\end{eqnarray}
for all $n$ large enough.
\end{lemma}
\begin{pf}
The uniform Harnack inequality implies that $q_j(x) \leq2C \rho(r)$
where $C = C(R/r)$ is the constant from the statement of part
\hyperlink{assumpharnack}{(2)} of Assumption \ref
{assumpfarawayharnack}. Let $F_j =
\{\tau_j(x) - \sigma_{j-1}(x) \leq T_\beta^U\}$. Arguing as in the
previous lemma and invoking uniform local transience, there exists
i.i.d. random variables~$\wt{F}_j(x)$ with $\p[\wt{F}_j(x) = 1] =
\delta= 1-\p[\wt{F}_j(x) = 0]$ that stochastically dominate $(\one
_{F_j(x)} : j)$ provided $R$ is sufficiently large. We let $\iota(j)$
be the $j$th smallest index $i$ such that $F_i(x)$ occurs. The lemma
now follows from an argument similar to that of Lemma \ref
{lemprobsuccessconcen}. Indeed, we can stochastically dominate~%
$q_{\iota(j)}(x)$ from below by i.i.d. random variables $L_j$ with $\E
L_j \geq(1-\delta) \ol{p}_{r,R}(x)$ and $L_j \leq10C \rho(r)$. By
Cram\'{e}r's theorem,
\[
\p\Biggl[ \prod_{j=1}^{k} (1-L_j) \geq\bigl(1-(1+\delta) \ol
{p}_{r,R}(x)\bigr)^{k}\Biggr]
\leq C\exp\bigl( - C \delta^2 \ol{p}_{r,R}(x) k / \rho(r) \bigr).
\]
The lemma now follows since, again by Cram\'{e}r's theorem,
\[
\p\bigl[ \iota\bigl((1-\delta)k\bigr) \geq k\bigr] \leq C \exp(-C \delta^2 k).
\]
\upqed\end{pf}

\subsection{\texorpdfstring{Proof of Theorem \protect\ref{thmthreshold}}{Proof of Theorem 1.3}}

We begin\vspace*{1pt} by showing that the points not visited by $X$ by time $\frac
{1}{2} T_\cov(G_n)$ are typically well separated, which in turn will
be helpful when we estimate the exponential moment in Proposition~\ref
{proptvbound}. To this end, we let $\delta> 0$ be arbitrary and assume
that $R > r,n_0,\varepsilon$ have been chosen so that for all $n \geq
n_0$ we have
\[
1-\delta\leq\frac{T_\cov(G_n)}{C_n^\varepsilon} \leq1+\delta.
\]
We may assume without loss of generality that $d_k^\varepsilon> 0$ for
all relevant $k$ and, in particular, that $|H_{n,k}^\varepsilon
|^{-\delta
} \to0$ for every $k$. Indeed, Lemmas \ref{lemhittingestimate} and~\ref{lemoccestimate} imply that $T_\hit(G_n) = \Theta(|V_n|)$,
consequently if ${\log}| H_{n,k}^\varepsilon| \to0$ as $n \to\infty$
then $T_\cov(H_{n,k}^\varepsilon)$ is negligible in comparison
to\vspace*{1pt}
$T_\cov(G_n)$. If $(G_n)$ satisfies Assumption \ref
{assumpfarawayharnack}\hyperlink{assumpfaraway}{(1)} we take $R_n^\gamma
$ as given there. Otherwise, we take $R_n^\gamma= \max\{ R > 0
\dvtx\break
{\max
_{x \in V_n}} |B(x,R)| \leq|V_n|^{\gamma}\}$.
\begin{lemma}
\label{lemwellseparated}
Let $\CR(t)$ denote the range of random walk at time $t$ and $\CL(t)
= V \setminus\CR(t)$. Letting
\[
M = \cases{
\displaystyle 20 \Delta_0 \sup_{n} \ol{\Delta}{}^R(G_n)/ (\delta\varepsilon
d^\varepsilon), &\quad if $\sup_{n} \ol{\Delta}(G_n) <
\infty$,\vspace*{2pt}\cr
\displaystyle 20 \Delta_0 / (\delta\varepsilon d^\varepsilon), &\quad otherwise,}
\]
and
\[
\CT_0 = \min\Bigl\{T \geq0 \dvtx\max_{x} |\CL(t) \cap B(x,R_n^\gamma
)| \leq M \Bigr\},
\]
we have\vspace*{1pt} that $\p[ \CT_0 > \frac{1+5\delta}{2} T_\cov(G_n)] =
o(1)$ provided $\gamma$ is sufficiently small,~$R$ is so large
that\vadjust{\goodbreak}
$\delta_{R,m} \leq1$, $d^\varepsilon= \min\{ d_k^\varepsilon\dvtx
d_k^\varepsilon> 0\}$ and $m =20/d^\varepsilon$. Furthermore, letting
\[
\CT_1 = \min\{T \geq0 \dvtx|\CL(t) \cap H_{n,k}^\varepsilon| \leq
|H_{n,k}^\varepsilon|^{1/2-\delta} \mbox{ for all } k \}
\]
we have that $\p[ \CT_1 > \frac{1+5\delta}{2} T_\cov(G_n)] = o(1)$.
\end{lemma}
\begin{pf}
First, suppose that $(G_n)$ has uniformly bounded maximal degree. Fix
$R > r$ and let $A$ be an $R$-net of $H_{n,k}^\varepsilon$. Fix $x \in
H_{n,k}^\varepsilon$ and suppose that $x_1,\ldots,x_\ell\in
B(x,R_n^\gamma) \cap H_{n,k}^\varepsilon\cap A$ are distinct. Lemma \ref
{lemcorrdecayupperbound} gives us
\[
\p\bigl[ x_1,\ldots, x_\ell\in\CL\bigl( (1+\delta)/2;G_n\bigr)\bigr] \leq(1+\delta
_{R,\ell}) |V_n|^{-(1+\delta)\ell d_k^\varepsilon/2 + \delta_{R,\ell}}.
\]
Consequently, a union bound yields
\begin{eqnarray*}
&&\p\bigl[ \bigl|\CL\bigl( (1+\delta)/2;G_n\bigr) \cap B(x,R_n^\gamma) \cap A\bigr| \geq\ell
\bigr]\\
&&\qquad\leq (1+\delta_{R,\ell}) |B(x,R_n^\gamma)|^\ell|V_n|^{-(1+\delta
)\ell d_k^\varepsilon/2 + \delta_{R,\ell}}\\
&&\qquad\leq (1+\delta_{R,\ell}) |V_n|^{(\gamma- (1+\delta)d_k^\varepsilon
/2)\ell+ \delta_{R,\ell}}.
\end{eqnarray*}
Hence, choosing $\gamma\leq d^\varepsilon/4$ the above is
$O(|V_n|^{-3})$. Since the number of disjoint $R$-nets necessary to
cover $H_{n,k}^\varepsilon$ is at most $\ol{\Delta}{}^R(G_n)$, the result
now follows from a union bound. In the case of unbounded maximal
degree, we can skip the step of subdividing the $H_{n,k}^\varepsilon$
into $R$-nets since in this case $\delta_{1,m} \to0$, otherwise the
proof is the same. The second claim is immediate from Markov's
inequality and Lem\-ma~\ref{lemcorrdecayupperbound}.
\end{pf}

We can now complete the proof of Theorem \ref{thmthreshold}. We will
handle the two cases depending on whether $(G_n)$ satisfies part
\hyperlink{assumpfaraway}{(1)}
or \hyperlink{assumpharnack}{(2)} of Assumption \ref
{assumpfarawayharnack}. Throughout, we let $N(x,T)$ be the number of
such excursions from $\partial B(x,r)$ to $\partial B(x,R)$ that have
occurred by time $T$.
\begin{pf*}{Proof of Theorem \ref{thmthreshold}, under Assumption \ref
{assumpfarawayharnack}\hyperlink{assumpharnack}{(2)}}
Let
\[
\CT_2 = \min\Biggl\{ T \geq0\dvtx\max_{x \in H_{n,k}^\varepsilon} \prod
_{k=1}^{N(x,T)} \bigl(1-q_j(x)\bigr) \leq|H_{n,k}^\varepsilon|^{-1/2-\delta}
\mbox{ for all } k \Biggr\}
\]
and set
%
%
\begin{equation}
\label{eqnstoppingtime}
\CT= \CT_0 \vee\CT_1 \vee\CT_2 \vee\biggl( \frac{1+5\delta
}{2}\biggr) T_\cov(G_n).
\end{equation}
Let $k_0(n)$ be a sequence so that $\liminf_{n \to\infty}
d_{k_0(n)}^\varepsilon\geq\delta_0 > 0$. For $x \in
H_{n,k_0(n)}^\varepsilon$, we have
\[
\biggl( \frac{1+3\delta}{2} \biggr) C_{n,k_0(n)}^\varepsilon\geq
\biggl(\frac{1+3\delta+O(\varepsilon)}{2}\biggr) \frac{{\delta_0 T_{r,R}(x)
\log}|V_n|}{4\rho(r)}
\]
for all $n$ large enough. Thus letting $M_{n,k_0(n)}^\varepsilon(x)
=(1+3\delta)/2 \cdot C_n^\varepsilon(x) / T_{r,R}(x)$, we have
\[
M_{n,k_0(n)}^\varepsilon(x) \geq\biggl( \frac{1+3\delta+ O(\varepsilon
)}{2} \biggr) \frac{{\delta_0 \log}|V_n|}{4 \rho(r)}.
\]
Now,
\begin{eqnarray*}
&&\p\bigl[ (1-\delta) T_{r,R}(x) M_{n,k_0(n)}^\varepsilon(x) \leq\tau
_{M_{n,k_0(n)}^\varepsilon}(x) \leq(1+\delta) T_{r,R}(x)
M_{n,k_0(n)}^\varepsilon(x)\bigr]\\
&&\qquad\geq 1- C\exp\biggl(- {\frac{C \delta_0 \delta^2}{\rho(r)}\log}
|V_n|\biggr)
\geq1 - O(|V_n|^{-100}),
\end{eqnarray*}
provided we choose $r$ large enough. Choosing $R > r$ sufficiently
large, Lem\-ma~\ref{lemprobsuccessbound} gives us
\[
\p\Biggl[ \prod_{j=1}^{M_{n,k_0(n)}^\varepsilon(x)}\bigl(1-q_j(x)\bigr) \geq
\bigl|H_{n,k_0(n)}^\varepsilon\bigr|^{-1/2-\delta}\Biggr]
\leq O(|V|^{-100}).
\]
Combining everything,
%
%
\begin{equation}
\label{eqnstoppingtimebound}
\p\biggl[ \CT\neq\biggl(\frac{1+5\delta}{2}\biggr) T_\cov
(G_n)\biggr] = o(1) \qquad\mbox{as } n \to\infty.
\end{equation}
Let $\mu$ be the probability on $\CX(G_n)$ given by first sampling
$\CR\subseteq V_n$ according to $\mu_0$, the measure on subsets of
$V_n$ given by running $X$ to time $(1+5\delta)/2 \cdot T_\cov(G_n)$,
then sampling $f|_{\CR}$ by marking with i.i.d. fair coins and
$f|_{V_n \setminus\CR} \equiv0$. Define $\wt{\mu}$ similarly
except by sampling $\CR\subseteq V_n$ according to $\wt{\mu}_0$, the
measure given by running $X$ up to time $\CT$ rather than $(1+5\delta
)/2 \cdot T_\cov(G_n)$. As a~consequence of~(\ref{eqnstoppingtimebound}),
\[
\|\mu- \wt{\mu}\|_{\mathrm{TV}} \leq\p\biggl[ \CT\neq\biggl(\frac
{1+5\delta}{2}\biggr) T_\cov(G_n)\biggr] = o(1) \qquad\mbox{as } n \to
\infty.
\]

Suppose we have two independent random walks $X,X'$ on $G_n$, each with
stationary initial distribution, and let $\CT,\CT'$ be stopping times
for each as in (\ref{eqnstoppingtime}). Let $\CR,\CR'$ be their
ranges at time $\CT,\CT'$, respectively, and $\CL= V_n \setminus\CR
$, $\CL' = V_n \setminus\CR'$. Let $q_j'(x)$ be the quantity
analogous to $q_j(x)$ for $X'$ and $\CG= \sigma(q_j'(x) \dvtx j \geq1)$.
The previous lemma implies that we can divide~$\CL$ into $M$ disjoint
sets $A_1,\ldots,A_M$ such that if $x,y \in A_\ell$ with $x \neq y$
then $d(x,y) \geq\break R_n^\gamma> R$. Consequently, letting $\CG(A_\ell)
= \otimes_{x \in A_\ell} \CG(x)$ we have
\begin{eqnarray*}
\E[\exp(\zeta|\CL\cap\CL' \cap A_\ell|) | \CG(A_\ell)]
&\leq& \prod_{x \in A_\ell} \Biggl(1 + e^{\zeta}\Biggl(\prod
_{j=1}^{N(x,\CT')} \bigl(1-q_j'(x)\bigr)\Biggr)\Biggr)\\
&\leq& \exp\biggl(e^{\zeta} \sum_k |H_{n,k}^\varepsilon|^{-\delta
}\biggr).
\end{eqnarray*}
Since $A_1,\ldots,A_M$ cover $\CL$, it follows from H\"older's
inequality that
%
%
\begin{eqnarray}\label{eqnholderbound}
\hspace*{68pt}\E\exp(\zeta|\CL\cap\CL'|)
&\leq&\biggl[ \exp\biggl(e^{\zeta M} \sum_k |H_{n,k}^\varepsilon|^{-\delta}\biggr)
\biggr]^{1/M} \nonumber\\[-8pt]\\[-8pt]
\hspace*{68pt}&\leq& 1 + 2\frac{\exp(\zeta M)}{M} \sum_k |H_{n,k}^{\varepsilon
}|^{-\delta} \nonumber.\hspace*{68pt}\qed
\end{eqnarray}
\noqed\end{pf*}
\begin{pf*}{Proof of Theorem \ref{thmthreshold}, under Assumption \ref
{assumpfarawayharnack}\hyperlink{assumpfaraway}{(1)}}
Let
\[
\CT_2 = \min\biggl\{T \geq0 \dvtx\max_k \max_{x \in H_{n,k}^\varepsilon
} \frac{{(1+2\delta)\log}|H_{n,k}^\varepsilon|}{2 N(x,T) \ol
{p}_{r,R}(x)} \leq1 \biggr\}
\]
and
%
%
\begin{equation}
\label{eqnstoppingtimesep}
\CT= \CT_0 \vee\CT_1 \vee\CT_2 \vee\biggl( \frac{1+5\delta
}{2}\biggr) T_\cov(G_n).
\end{equation}
It follows from Lemmas \ref{lemexcursionconcen} and \ref
{lemcovertime} and the definition of $H_{n,k}^\varepsilon$ that
%
%
\begin{equation}
\label{eqnstoppingtimeboundsep}
\p\biggl[ \CT\neq\biggl(\frac{1+5\delta}{2}\biggr) T_\cov
(G_n)\biggr] = o(1) \qquad\mbox{as } n \to\infty.
\end{equation}
Let $\mu$ be the probability on $\CX(G_n)$ given by first sampling
$\CR\subseteq V_n$ according to $\mu_0$, the measure on subsets of
$V_n$ given by running $X$ to time $(1+5\delta)/2 \cdot T_\cov(G_n)$,
then sampling $f|_\CR$ by marking with i.i.d. fair coins and $f|_{V_n
\setminus\CR} \equiv0$. Define $\wt{\mu}$ similarly except by
sampling $\CR\subseteq V_n$ according to $\wt{\mu}_0$, the measure
given by running $X$ up to time $\CT$ rather than $(1+5\delta)/2
\cdot T_\cov(G_n)$. As a~consequence of~(\ref{eqnstoppingtimeboundsep}),
\[
\|\mu- \wt{\mu}\|_{\mathrm{TV}} \leq\p\biggl[ \CT\neq\biggl(\frac
{1+5\delta}{2}\biggr) T_\cov(G_n)\biggr] = o(1) \qquad\mbox{as } n \to
\infty.
\]

Suppose we have two independent random walks $X,X'$ on $G_n$, each with
stationary initial distribution, and let $\CT,\CT'$ be stopping times
for each as in~(\ref{eqnstoppingtimesep}). Using the same notation as
the previous proof, by the definition of~$\CT_2'$, we have
%
%
\begin{eqnarray}\label{eqnexpboundsep}
&&\E\bigl[ \E[\exp(\zeta|\CL\cap\CL' \cap A_\ell|) | \CG(A_\ell
)] \bigr]
\nonumber\\
&&\qquad\leq \E\prod_{x \in A_\ell} \Biggl(1 + e^{\zeta}\Biggl(\prod
_{j=1}^{N(x,\CT')} \bigl(1-q_j'(x)\bigr)\Biggr)\Biggr) \\
&&\qquad\leq \E\prod_{x \in A_\ell} \Biggl(1 + e^{\zeta}\Biggl(\prod
_{j=1}^{N(x)} \bigl(1-q_j'(x)\bigr)\Biggr)\Biggr),\nonumber
\end{eqnarray}
where $N(x) = {(1+2\delta) \log}|H_{n,k}^\varepsilon| / 2 \ol
{p}_{r,R}(x)$ and $k$ is such that $x \in H_{n,k}^\varepsilon$. Let
\[
\wt{N}(x) = (1-\delta) N(x) \geq\frac{{(1+\delta/2) \log}
|H_{n,k}^\varepsilon|}{2 \ol{p}_{r,R}(x)}.
\]
Observe that (\ref{eqnexpboundsep}) is bounded by
\[
\E\prod_{x \in A_\ell} \Biggl(1 + e^{\zeta}\Biggl(\prod_{j=1}^{\wt
{N}(x)} \bigl(1-q_{i(j)}'(x)\bigr) + \one_{\{I(\wt{N}(x)) > N(x)\}}
\Biggr)\Biggr).
\]
As $A_\ell$ satisfies the hypotheses of Lemma \ref
{lemiidreplacement}, this is in turn bounded by
\begin{eqnarray*}
&& \E\prod_{x \in A_\ell} \Biggl(1 + e^{\zeta}\Biggl(\prod
_{j=1}^{\wt{N}(x)} \bigl(1-(1-\delta/4)Q_j'(x)\bigr)\Biggr) +
O(|V_n|^{-100})\Biggr)\\
&&\qquad\leq \exp\biggl(e^{\zeta} \sum_k |H_{n,k}^\varepsilon|^{-\delta
}\biggr).
\end{eqnarray*}
The theorem now follows from H\"older's inequality, as in the previous proof.
\end{pf*}

\subsection{The lamplighter}

\mbox{}

\begin{pf*}{Proof of Theorem \ref{thmlamplighter}}
This is proved by making several small modifications to the proof of
Theorem \ref{thmthreshold}. Namely, rather than considering the range
of $X$ run up to time $\CT$ as in either (\ref{eqnstoppingtime}) or
(\ref{eqnstoppingtimesep}), one considers the range $\wt{\CR}(x)$ of
$X$ run up to time $\CT$, conditioned on the event $\{X(\CT) = x\}$
for a given point $x$. Exactly the same argument shows that the total
variation distance of the law $\wt{\mu}_x$ on markings $\CX(G_n)$
induced by i.i.d. coin flips on $\wt{\CR}(x)$ and $0$ on $(\wt{\CR
}(x))^c$ from the uniform measure on $\CX(G_n)$ is $o(1)$. This
implies that the law $\mu_x$ on markings of $\CX(G_n)$ given by
i.i.d. coin flips on the range $\CR(x)$ of $X$ run up to time $T =
\frac{1+\varepsilon}{2} T_\cov(G_n)$, conditioned on $\{X(T)=x\}$, and
the uniform measure is $o(1)$. At time $T$, the random walk is well
mixed, from which the result is clear.
\end{pf*}

\section{Further questions}

\begin{longlist}[2.]
\item[1.] Theorem \ref{thmthreshold} yields a wide class of examples where
the threshold for indistinguishability is at $\frac{1}{2} T_\cov$,
and $\Z_n^2$ is an example where the threshold is at $T_\cov$. Does
there exist a sequence $(G_n)$ of vertex transitive graphs where the
threshold is at $\alpha T_\cov(G_n)$ for $\alpha\in(1/2,1)$?
\item[2.] Our statistical test for uniformity is only valid for $\alpha>
1/2$. For $\alpha\leq1/2$, the natural reference measure is i.i.d.
markings conditioned on the number of zeros being on the order of
$|V|^{1-\alpha}$. Can analogous results be proved in this setting?
\item[3.] Our definition of uniform local transience is given in terms of
 Green's function summed up to the \textit{uniform mixing time}. Does
it suffice to assume only the uniform decay of
\[
g(x,y;G) = \sum_{t=1}^T p^t(x,y;G),
\]
where $T = T_\mix(G)$ or even $T = T_{{\mathrm{rel}}}(G)$?
\item[4.] The complete graph $K_n$ does not satisfy the hypotheses of
Theorem~\ref{thmthreshold} yet the lamplighter walk on $K_n$ has a
threshold at $\frac{1}{2} T_\cov(K_n)$. Is there a more general
theorem allowing for a unified treatment of this case?
\end{longlist}

\section*{Acknowledgment}
J. Miller thanks the Theory Group at Microsoft Research for
support through a summer internship, during which the research for this article
was conducted.


%

%
\printaddresses

\end{document}